\documentclass[12pt,francais]{article}
\usepackage{amsmath}\usepackage{epsf,amsfonts,amsthm}\usepackage{amscd,amssymb}
\usepackage{xcolor,epic,eepic}\usepackage{epsfig}\usepackage{faktor}
\usepackage{fontenc,indentfirst, delarray,amsfonts,amsmath,amssymb}
\usepackage{rotating}
\usepackage{mathdots}
\usepackage[T1]{fontenc}
\usepackage[matrix,arrow,curve]{xy}
\usepackage{amsmath}
\usepackage{amssymb}
\usepackage{amsthm}
\usepackage{amscd}
\usepackage{amsfonts}
\usepackage{graphicx}
\usepackage{fancyhdr}
\usepackage{dsfont,texdraw}
\usepackage{amsmath}\usepackage{epsf,amsfonts,amsthm}\usepackage{amscd,amssymb}
\usepackage{xcolor,epic,eepic}\usepackage{epsfig}
\usepackage{fontenc,indentfirst, delarray,amsfonts,amsmath,amssymb}
\usepackage{rotating}

\usepackage{amscd}
\usepackage{amsmath}
\usepackage{amsthm}
\usepackage{mathrsfs}
\usepackage{latexsym}
\usepackage{amssymb}
\usepackage{amsfonts}
\theoremstyle{plain}

\usepackage{tikz}
\usetikzlibrary{arrows,chains,matrix,positioning,scopes,snakes}
\makeatletter
\tikzset{join/.code=\tikzset{after node path={%
\ifx\tikzchainprevious\pgfutil@empty\else(\tikzchainprevious)%
edge[every join]#1(\tikzchaincurrent)\fi}}}
\makeatother

\tikzset{>=stealth',every on chain/.append style={join},
         every join/.style={->}}
\tikzset{
    >=stealth',
    punkt/.style={
           rectangle,
           rounded corners,
           draw=black, very thick,
           text width=6.5em,
           minimum height=2em,
           text centered},
    pil/.style={
           ->,
           thick,
           shorten <=2pt,
           shorten >=2pt,}
}

\usepackage[T1]{fontenc}
\newcommand{\bee}{\begin{enumerate}}
\newcommand{\eee}{\end{enumerate}}
\newcommand{\benn}{\begin{equation*}}
\newcommand{\eenn}{\end{equation*}}
\newcommand{\be}{\begin{equation}}
\newcommand{\ee}{\end{equation}}
\newcommand{\bean}{\begin{eqnarray}}
\newcommand{\eean}{\end{eqnarray}}
\newcommand{\bea}{\begin{eqnarray*}}
\newcommand{\eea}{\end{eqnarray*}}

\newcommand{\N}{\mathbb{N}}
\newcommand{\Z}{\mathbb{Z}}

\newcommand{\cF}{{\cal F}}
\newcommand{\n}{\nabla}

\newcommand{\op}[1]{\!\!\mathop{\rm ~#1}\nolimits}

\newcommand{\id}{\op{id}}

\mathchardef\za="710B  
\mathchardef\zb="710C  
\mathchardef\zg="710D  
\mathchardef\zd="710E  
\mathchardef\zve="710F 
\mathchardef\zz="7110  
\mathchardef\zh="7111  
\mathchardef\zy="7112 

\mathchardef\zi="7113  
\mathchardef\zk="7114  
\mathchardef\zl="7115  
\mathchardef\zm="7116  
\mathchardef\zn="7117  
\mathchardef\zx="7118  
\mathchardef\zp="7119  
\mathchardef\zr="711A  
\mathchardef\zs="711B  
\mathchardef\zt="711C  
\mathchardef\zu="711D  
\mathchardef\zf="711E 
\mathchardef\zq="711F  
\mathchardef\zc="7120  
\mathchardef\zw="7121  
\mathchardef\ze="7122  
\mathchardef\zvy="7123  
\mathchardef\zvw="7124  
\mathchardef\zvr="7125 
\mathchardef\zvs="7126 
\mathchardef\zvf="7127  
\mathchardef\zG="7000  
\mathchardef\zD="7001  
\mathchardef\zY="7002  
\mathchardef\zL="7003  
\mathchardef\zX="7004  
\mathchardef\zP="7005  
\mathchardef\zS="7006  
\mathchardef\zU="7007  
\mathchardef\zF="7008  
\mathchardef\zW="700A  

\newcommand{\cyclic}{\mathop{\kern0.9ex{{+}
\kern-2.15ex\raise-.25ex\hbox{\Large\hbox{$\circlearrowright$}}}}\limits}

 \newcommand{\cA}{{\cal A}}
 
 \newcommand{\cD}{{\cal D}}
 \newcommand{\cO}{{\cal O}}
 \newcommand{\cB}{{\cal B}}

 \newcommand{\cI}{{\cal I}}

\newtheorem{rem}{Remark}
\newtheorem{theo}{Theorem}
\newtheorem{prop}{Proposition}
\newtheorem{lem}{Lemma}
\newtheorem{cor}{Corollary}

\newtheorem{defi}{Definition}

\newtheorem{lemma}{Lemma}

\newcommand{\h}{\op{Hom}}
\newcommand{\0}{\otimes}

\DeclareMathAlphabet{\mathpzc}{OT1}{pzc}{m}{it}

\usepackage{stackrel}
\usepackage{epigraph}

\usepackage{lmodern}
\usepackage{here}
\usepackage{url}
\usepackage{hyperref}

\usepackage{hyperref}
\usepackage{tikz-cd}

\def\lim{\mathop{\rm lim}\nolimits}

\newcommand{\xdownarrow}[1]{%
  {\left\downarrow\vbox to #1{}\right.\kern-\nulldelimiterspace}
}
\usepackage{eucal}

\def\clap#1{\hbox to 0pt{\hss#1\hss}}

\begin{document}

\title{\Large{{\it The Tor Spectral Sequence and Flat Morphisms\\ in Homotopical $\cD$-Geometry}}}
\author{Alisa Govzmann, Damjan Pi\v{s}talo, and Norbert Poncin}
\maketitle

\begin{abstract} \footnotesize{Homotopical algebraic $\cD$-geometry combines aspects of homotopical algebraic geometry \cite{TV05,TV08} and $\cD$-geometry \cite{BD04}. It was introduced in \cite{HAC} as a suitable framework for a coordinate-free study of the Batalin-Vilkovisky complex and more generally for the study of non-linear partial differential equations and their symmetries \cite{CG1,P11}. In order to consolidate the foundation of the theory, we have to prove that the standard methods of linear and commutative algebra are available in the context of homotopical algebraic $\cD$-geometry, and we must show that in this context the étale topology is a kind of homotopical Grothendieck topology and that the notion of smooth morphism is, roughly speaking, local for the étale topology. The first half of this work was done in \cite{HAC}. The remaining part covers the study of étale and flat morphisms in the category of differential graded $\cD$-algebras and is based on the $\op{Tor}$ spectral sequence which connects the graded Tor functors in homology with the homology of the derived tensor product of two differential graded $\cD$-modules over a differential graded $\cD$-algebra.}\end{abstract}

\small{\vspace{2mm} \noindent {\bf MSC 2020}: 18G15, 18E35, 18N40, 14A30 \medskip

}

\setcounter{tocdepth}{1}
\tableofcontents

\thispagestyle{empty}

\section{Introduction}\label{Intro}

Although deep insights into theoretical physics are often gained from calculations in coordinates, since the discovery of the general theory of relativity, geometers have been striving again for coordinate-free concepts and results. As many mathematical models reduce real-world problems to solving systems of partial differential equations (PDE-s), a coordinate-independent theory of PDE-s and related aspects is essential. In \cite{Vino}, A. Vinogradov and his co-workers propose a cohomological analysis of PDE-s, while A. Beilinson and V. Drinfeld advocate in \cite{BD04} working in the algebraic $\cD$-geometric setting, where $\cD$ is the non-commutative ring of total sections of the sheaf of differential operators of a smooth affine algebraic variety \cite{HAC,Gil06}. The solutions of a system of non-linear PDE-s are related to diffieties in the first approach and to $\cD$-schemes in the second. To treat the moduli space of solutions modulo symmetries, we have to replace $\cD$-schemes with $\cD$-stacks and more generally with derived $\cD$-stacks. In the functor-of-points approach to geometry \cite{TVHD,CCF,Schwarz}, derived $\cD$-stacks are those functors or presheaves $F:{\tt DG\cD A}\to{\tt SSet}$ from the model category of differential non-negatively graded $\cD$-algebras into the model category of simplicial sets, which are also sheaves. This suggests combining the $\cD$-geometry of Beilinson and Drinfeld with the homotopical algebraic geometry in the sense of B. Toën and G. Vezzosi \cite{TV05}, \cite{TV08}. In a series of papers \cite{KTRCR}, \cite{HAC}, \cite{PP} and \cite{HomotopyfiberSeq}, the authors present homotopical algebraic $\cD$-geometry as a suitable framework for the study of PDE-s and their symmetries. In fact, this new geometry provides in particular a convenient method for encoding total derivatives and leads to a covariant description of the classical Batalin-Vilkovisky complex, which appears as a special case of general constructions. Further evidence for this point of view can be found in \cite{P11}, \cite{Paugam1}.\medskip

In order to implement these ideas and in particular to define the aforementioned sheaf condition, we have to prove that the tuple $({\tt DG\cD M, DG\cD M, DG\cD A},\zt,\mathbf{P})$ is a homotopical algebraic geometric context (HAG context) in the sense of \cite{TV08}. Here $\tt DG\cD M$ is the symmetric monoidal model category of differential non-negatively graded $\cD$-modules (for $\cD$-modules, see for instance \cite{Schn}), $\zt$ is a suitable model pre-topology on the opposite category $\tt DG\cD A^{\op{op}}$ and $\mathbf{P}$ is a class of morphisms of $\tt DG\cD A^{\op{op}}$ that is compatible with $\zt$ (we can think of $\zt$ as the étale topology and of $\mathbf{P}$ as the class of smooth morphisms). In \cite{HAC} the first part of the HAG context conditions was proved, i.e., we showed that the triplet $({\tt DG\cD M, DG\cD M, DG\cD A})$ is a homotopical algebraic context (HA context) \cite{TV08}, which guarantees that essential tools from linear and commutative algebra are still available. To prove the remaining part, we need to extend Quillen's Tor spectral sequence \cite{Quill} and the correct homotopy invariant notion of flatness \cite{TV08} to the $\cD$-geometric context. The generalization of the Tor spectral sequence connects the homology of the derived tensor product of two modules $M,N\in \tt DG\cD M$ over an algebra $\cA\in \tt DG\cD A$ with the Tor of the corresponding homologies, which appears as the second sheet of a first quadrant spectral sequence in $\tt \cD M\,$. If $M$ is flat in the just mentioned homotopy invariant sense and more precisely in the sense that the derived tensor product with $M$ preserves homotopy pullbacks, the Tor spectral sequence collapses at its second sheet, giving rise to the isomorphism of $H\cA$-modules $$H M\otimes_{H\cA}H N\cong H(M\otimes^\mathbb{L}_\cA N)\;.$$ It is also worth mentioning that, in contrast to the present work, the corresponding parts in \cite{Quill} and \cite{TV08} are written in the context of simplicial modules. Besides enriching the theory by including the action of differential operators, our contribution translates Quillen's spectral sequence and flatness from the world of simplicial modules into the language of differential graded modules. While simplicial modules are indispensable when working over a base ring of possibly positive characteristic, in derived algebraic geometry one usually prefers differential graded modules when the base ring has characteristic zero. Although the general philosophy of the proofs is the same in both situations, the technical reasoning is different.\medskip

The paper is organized as follows.\medskip

Section \ref{Zoom} establishes the conventions used throughout this text, particularly in relation to the definition of model categories, their homotopy category, and the definition of derived functors, as well as the replacements used to compute them. Theorem \ref{Fundamental0} makes it easier to navigate the thicket of similar but different concepts, definitions, and results -- see for instance \cite{DS,Hir,Ho99,JL}. Indeed, although the facts presented in this theorem are by no means new, to our knowledge it is difficult to find a similar result in the literature.\medskip

In Section \ref{DTPF} we prove that in our category ${\tt Mod}(\cA)$ of modules $M,N$ in $\tt DG\cD M$ over $\cA\in\tt DG\cD A$ (i.e., of modules in the category of differential non-negatively graded $\cD$-modules over a commutative monoid in this symmetric monoidal category), in order to compute the derived tensor product $M\0^{\mathbb{L}}_\cA N\,,$ we can choose a cofibrant replacement of $M$ and $N,$ of $M$ only, or $N$ only, just as in classical homological algebra.\medskip

The family of Tor functors is usually defined in the setting of modules over a ring. Considering what we said above, this family appears here in the context of (homology) modules in ${\tt G\cD M}$ over $\mathfrak{A}\in\tt G\cD A\,.$ Hence the value on objects of each one of the Tor functors is a graded $\cD$-module over the graded $\cD$-algebra $\mathfrak{A}\,.$ The modules of this category ${\tt Mod}(\mathfrak{A})$ are {\it not} graded modules over the single ring of differential operators with coefficients in $\mathfrak{A}\,.$ Instead, we are dealing in addition to the $\mathfrak{A}$-action with a compatible action by the non-commutative ring $\cD$ of differential operators. Further, the $\cD$-context is a rather special one, among other things because the tensor product $M\0_\cA N$ (for notations, see preceding paragraph) is a quotient of the tensor product $M\0 N$ of differential graded $\cD$-modules, which however is not the tensor product over $\cD\,$, but that over the commutative ring $\cO$ of functions. In Section \ref{GraTor} we show that the just mentioned graded Tor functors do exist, which includes checking that the categories ${\tt Mod}(\cA)$ and ${\tt Mod}(\mathfrak{A})$ are abelian, that the latter has enough projectives, and that the partial tensor product functor over it is right exact.\medskip

In Section \ref{TSSDGeo} we give the $\cD$-geometric generalization of Quillen's Tor spectral sequence -- a particular Künneth spectral sequence used by Quillen within the framework of simplicial modules over a simplicial ring. We explain why Quillen's original proof can be generalized to our setting. In particular, we use Sullivan differential graded $\cD$-modules over differential graded $\cD$-algebras \cite{KTRCR,HAC} and graded $\cD$-modules over graded $\cD$-algebras which are free over a graded set. In addition, we prove a strong Künneth theorem in this context. \medskip

In Section \ref{fm1} we give all the details needed for a rigorous homotopy invariant definition of flatness (see above). More precisely, on the one hand, the homotopy pullback in $\tt Mod(\cA)$ is a functor from the homotopy category ${\tt Ho(CoSpan(Mod(\cA)))}$ of cospan diagrams in $\cA$-modules to the category $\tt Ho(Mod(\cA))$. On the other hand, the derived tensor product with $M\in \tt Mod(\cA)$ is an endofunctor of $\tt Ho(Mod(\cA))$ and induces an endofunctor of $\tt CoSpan( Ho(Mod(\cA) ))$, but not an endofunctor of the non-equivalent homotopy category $\tt Ho(CoSpan(Mod(\cA)))$. Therefore, it is not clear at first what it means for the derived tensor product to preserve homotopy pullbacks. Also, similar to usual pullbacks, for a functor $\cF$ to preserve homotopy pullbacks it is not sufficient that the $\cF$-image of the homotopy pullback of a diagram is isomorphic to the homotopy pullback of the $\cF$-image, but this `isomorphism' must be realized by a specific universal weak equivalence. Section \ref{fm1} solves these difficulties and leads to a completely rigorous concept of flatness (Definition \ref{FlatDefi}). It is further shown (Corollary \ref{FlatSeq}) that for an $\cA$-module $M$ to be flat it is sufficient that the derived tensor product with $M$ preserves homotopy fiber sequences in the sense of \cite{HomotopyfiberSeq}. As with the preservation of homotopy pullbacks, a few subtleties have to be considered in order to give the preservation of homotopy fiber sequences a rigorous meaning (Definition \ref{SeqPul}).\medskip

In Section \ref{fm2} we prove that an $\cA$-module in the category of differential graded $\cD$-modules is flat if and only if it is strongly flat (for a similar result in the case of a commutative base ring, see \cite{TV08}).\medskip

We conclude with a brief outlook in Section \ref{outlook}.

\section{Zooming in on derivatives in model categories}\label{Zoom}

This section establishes notations and conventions and ensures that the present text is self-contained. Theorem \ref{Fundamental0} is a model categorical generalization of the well-known homological algebraic fact that the values of the left derived functors of a covariant right exact functor from an abelian category with enough projectives to an abelian category, are independent, up to canonical isomorphism, of the projective resolution used.\medskip

We assume that the reader is familiar with model categories. We adopt the definition of a model category that is used in \cite{Hir}. More precisely, a model category is a category $\tt M$ that is equipped with three classes of morphisms called weak equivalences, fibrations and cofibrations. The category $\tt M$ has all {\it small} limits and colimits and the 2-out-of-3 axiom, the retract axiom and the lifting axiom are satisfied. Moreover $\tt M$ admits a {\it functorial} cofibration - trivial fibration factorization system (Cof - TrivFib factorization) and a functorial trivial cofibration - fibration factorization system (TrivCof - Fib factorization).\medskip

Let now $(\za,\zb)$ be any functorial Cof - TrivFib factorization. For every $X\in\tt M\,,$ it factors the map $i_X:0\to X$ out of the initial object of $\tt M$ into a cofibration $\za(i_X)$ followed by a trivial fibration $q_X:=\zb(i_X)\,$: $$i_X:0\rightarrowtail QX\stackrel{\sim}{\twoheadrightarrow}X\;.$$ Regardless of the factorization $$i_X:0\rightarrowtail CX\stackrel{\sim}{\to}X$$ of $i_X$ into a cofibration followed by a weak equivalence $c_X$ considered, we refer to $CX$ as {\it a} {\it cofibrant replacement} of $X\,.$ The object $QX$ we call a {\it cofibrant F-replacement} of $X$ (or just a cofibrant replacement if we do not want to stress that $q_X$ is a fibration). From the fact that the factorization $(\za,\zb)$ is functorial it follows that $Q$ is an endofunctor of $\tt M\,$. Moreover $q_X:QX\to X$ is functorial in $X\,$: $q$ is a natural transformation $q:Q\Rightarrow\op{id}_{\tt M}$ from the {\it cofibrant replacement functor} $Q$ to the identity functor $\op{id}_{\tt M}$ \cite{Ho99}. Instead of the cofibrant F-replacement functor $Q$ that is {\bf global}ly defined by the functorial factorization $(\za,\zb)\,,$ we will also use {\bf local} / object-wise cofibrant replacements $CX$ or {\it local cofibrant F-replacements} $\tilde{C}X$ such that the map $c_X$ in the factorization $$i_X:0\rightarrowtail \tilde{C}X\stackrel{\sim}{\twoheadrightarrow}X$$ is $\op{id}_X$ if $X$ is already cofibrant \cite{Quill}. It is important to remember that if for every $X$ we choose such a local cofibrant F-replacement and if $f:X\to Y\,,$ there exists a lifting $\tilde{C}f:\tilde{C}X\to \tilde{C}Y\,:$

\begin{equation}\label{CTilde} \begin{tikzpicture}
 \matrix (m) [matrix of math nodes, row sep=3em, column sep=3em]
   {  \stackrel{}{0}  & & \tilde{C}Y  \\
      \tilde{C}X & \stackrel{}{X} & \stackrel{}{Y}  \\ };
 \path[->]
 (m-1-1) edge [>->] (m-1-3)
 (m-1-1) edge [>->] (m-2-1)
 (m-2-1) edge [->>] node[below] {\small{$\;\;{}_{\widetilde{}}\,\;c_X$}} (m-2-2)
 (m-1-3) edge [->>] node[auto] {\small{$\;{}_{\widetilde{}}\;\;c_Y$}}(m-2-3)
 (m-2-2) edge node[below] {\small{$f$}} (m-2-3)
 (m-2-1) edge [->, dashed] node[auto] {\small{$\tilde{C}f$}} (m-1-3);
\end{tikzpicture}
\end{equation}

By homotopy category of a model category $\tt M$ in this work we mean the Quillen homotopy category $\tt Ho(M)\,,$ which is the strong localization ${\tt M}[[W^{-1}]]$ of $\,\tt M$ at its class $W$ of weak equivalences \cite{Ho99}. Like any localization, the category $\tt Ho(M)$ comes along with a localization functor $L_{\tt M}:\tt M \to Ho(M)$ that sends weak equivalences to isomorphisms. By strong localization we mean that if $\tt C$ is another category and $F: \tt M \to C$ another functor that sends weak equivalences to isomorphisms, then there exists a {\it unique} functor ${\tt Ho}(F):{\tt Ho(M)}\to{\tt C}$ such that the resulting triangle commutes {\it on the nose}, i.e., $F={\tt Ho}(F)\circ L_{\tt M}\,$. The objects of ${\tt Ho(M)}$ are the objects of $\tt M$ and the morphisms of ${\tt Ho(M)}$ are the alternating finite strings $$[S]=[\to\,\rightsquigarrow\,\to\,\rightsquigarrow\,\cdots]\quad\text{and}\quad[S]=[\rightsquigarrow\,\to\,\rightsquigarrow\,\to\,\cdots]$$ of morphisms $f:X\to Y$ of $\tt M$ and formal reversals $w^{-1}:Y\rightsquigarrow Z$ of weak equivalences $w:Z\to Y$ of $\tt M\,.$ Here the class $[\cdot]$ refers to the identifications of the empty string $1_X$ at X, the concatenation string $f,g:X\to Y\to Z$ of composable $\tt M$-morphisms, the concatenation strings $w,w^{-1}:Z\to Y\rightsquigarrow Z$ and $w^{-1},w:Y\rightsquigarrow Z\to Y$ with the $\tt M$-morphisms $\id_X,$ $g\circ f$, $\id_Z$ and $\id_Y\,,$ respectively. The unique factorization ${\tt Ho}(F)$ of $F$ through $\tt Ho(M)$ is defined on objects $X\in\tt Ho(M)$ by ${\tt Ho}(F)(X)=F(X)$ and on $\tt Ho(M)$-morphisms $[S]$ by $${\tt Ho}(F)[f:X\to Y]=F(f:X\to Y)\quad\text{and}\quad {\tt Ho}(F)[w^{-1}:Y\rightsquigarrow Z]=(F(w:Z\to Y))^{-1}\;.$$

From what has been said in the previous paragraph it follows that a functor $F:{\tt M}\to {\tt N}$ between model categories $\tt M$ and $\tt N\,,$ which preserves {\it all} weak equivalences can be viewed canonically as a functor between the homotopy categories. Indeed, the composite $L_{\tt N}\circ F:{\tt M}\to {\tt Ho(N)}$ sends weak equivalences to isomorphisms and therefore factors uniquely through the homotopy category of $\tt M\,$, i.e., there is a unique functor \be\label{CanDerFun}{\tt Ho}(L_{\tt N}\circ F):{\tt Ho(M)}\to {\tt Ho(N)}\ee such that $${\tt Ho}(L_{\tt N}\circ F)\circ L_{\tt M}=L_{\tt N}\circ F\;.$$ However, requiring that functors between model categories respect the {\it entire} model structure is far too strong a requirement. Usually, functors $F$ between model categories are left (or right) Quillen functors, i.e., they preserve all cofibrations and all trivial cofibrations (resp., all fibrations and all trivial fibrations). Since a left Quillen functor $F$ thus sends trivial cofibrations between cofibrant objects to weak equivalences, it follows from Brown's lemma that it sends all {\it weak equivalences between cofibrant objects} to weak equivalences. Hence, if $Q_{\tt M}$ denotes a cofibrant replacement functor of $\,{\tt M}$, since due to the 2-out-of-3 axiom $Q_{\tt M}$ preserves weak equivalences, the functor $L_{\tt N}\circ F\circ Q_{\tt M}:{\tt M}\to {\tt Ho(N)}$ sends weak equivalences to isomorphisms and therefore factors uniquely through the homotopy category of $\tt M\,$, i.e., there is a unique functor \be\label{DerFun}{\tt Ho}(L_{\tt N}\circ F\circ Q_{\tt M}):{\tt Ho}({\tt M})\to {\tt Ho(N)}\ee such that $${\tt Ho}(L_{\tt N}\circ F\circ Q_{\tt M})\circ L_{\tt M} = L_{\tt N}\circ F\circ Q_{\tt M}\;.$$ We refer to this functor as `the' left derived functor of $F\,.$ Since its definition is based on the universal property of the strong localization $\tt Ho(M)$ and it is computed using the cofibrant replacement functor $Q_{\tt M}\,,$ we also speak of the {\bf strongly universal left derived functor} of $F$ and denote this functor $\mathbb{L}^{\op{S}}_Q F\,:$ \be\label{DefDerF} \mathbb{L}^{\op{S}}_Q F:= {\tt Ho}(L_{\tt N}\circ F\circ Q_{\tt M})\;.\ee Alternatively one can define the left derived functor of $F$ as the right Kan extension of $L_{\tt N}\circ F$ along $L_{\tt M}\,.$ We refer to this derived functor as the {\bf Kan extension left derived functor} of $F$ and denote it by $\mathbb{L}^{\op{K}} F\,.$ The right derived functors $\mathbb{R}^{\op{S}}_R G$ ($R=R_{\tt M}$: fibrant replacement functor) and $\mathbb{R}^{\op{K}} G$ of a right Quillen functor $G$ are defined dually.\medskip

We first clarify the relationship between \eqref{CanDerFun} and \eqref{DerFun}:

\begin{prop}
If $F\in\tt Fun(M,N)$ is a functor between model categories that preserves \emph{all} weak equivalences, the whiskering $\zy:=L_{\tt N}\star(F\star q_{\tt M})$ is a canonical natural isomorphism $$\zy:{\tt Ho}(L_{\tt N}\circ F\circ Q_{\tt M})\stackrel{\cong}{\Rightarrow}{\tt Ho}(L_{\tt N}\circ F)\;.$$ We refer to the latter by simply writing $$\mathbb{L}^{\op{S}}_QF:={\tt Ho}(L_{\tt N}\circ F\circ Q_{\tt M})\doteq{\tt Ho}(L_{\tt N}\circ F)\;.$$
\end{prop}

\begin{proof}
As $q_{\tt M}:Q_{\tt M}\stackrel{\sim}{\Rightarrow}\id_{\tt M}$ is a natural weak equivalence, i.e., a natural transformation that is objectwise a weak equivalence, the whiskering $\zy: = L_{\tt N}\star(F\star q_{\tt M})$ is a natural isomorphism $$\zy:{\tt Ho}(L_{\tt N}\circ F\circ Q_{\tt M})\circ L_{\tt M}\stackrel{\cong}{\Rightarrow}{\tt Ho}(L_{\tt N}\circ F)\circ L_{\tt M}\;,$$ because for every $X\in\tt M$ the morphism $\zy_X:=L_{\tt N}(F(q_{{\tt M},X}))$ is an isomorphism. That $\zy$ is a natural isomorphism $$\zy:{\tt Ho}(L_{\tt N}\circ F\circ Q_{\tt M})\stackrel{\cong}{\Rightarrow}{\tt Ho}(L_{\tt N}\circ F)\;$$ follows from \cite[Lemma 1]{CompTheo} which is a rather straightforward consequence of the above description of $\tt Ho(M)\,.$
\end{proof}

The next theorem addresses the question of stability of a derived functor with respect to a change of definition (Kan extension versus strongly universal) and with respect to a change of the type of cofibrant replacement used to compute it (local versus global):

\begin{theo}[\cite{Models}]\label{Fundamental0}
If $F\in\tt Fun(M,N)$ is a functor between model categories that sends weak equivalences between cofibrant objects to weak equivalences, its Kan extension left derived functor $$\mathbb{L}^{\op{K}}F\in\tt Fun(Ho(M),Ho(N))$$ and its strongly universal left derived functor $$\mathbb{L}^{\op{S}}F\in\tt Fun(Ho(M),Ho(N))$$ (\cite{CompTheo}) exist and we have
\be\label{Fundamental1}\mathbb{L}^{\op{K}}F\doteq {\tt Ho}(L_{\tt N}\circ F\circ\tilde{C}_{\tt M})\doteq\mathbb{L}_Q^{\op{S}}F:={\tt Ho}(L_{\tt N}\circ F\circ Q_{\tt M})\stackrel{\cong}{\Rightarrow} \mathbb{L}^{\op{S}}F\;,\ee where $\tilde{C}_{\tt M}$ is a local cofibrant F-replacement, $Q_{\tt M}$ is a cofibrant F-replacement functor and ${\tt Ho}$ the unique on-the-nose factorization through $\tt Ho(M)\,$. This implies that \be\label{Fundamental2}\mathbb{L}^{\op{K}}F\circ L_{\tt M}\doteq L_{\tt N}\circ F\circ\tilde{C}_{\tt M}\doteq\mathbb{L}_Q^{\op{S}}F\circ L_{\tt M}=L_{\tt N}\circ F\circ Q_{\tt M}\stackrel{\cong}{\Rightarrow}\mathbb{L}^{\op{S}}F\circ L_{\tt M}\;,\ee where $\doteq$ denotes a canonical natural isomorphism and $\stackrel{\cong}{\Rightarrow}$ a not necessarily canonical natural isomorphism.
\end{theo}

\noindent Hence, for every $X\in\tt M\,,$ the value of the derived functor at $L_{\tt M}X=X\in\tt Ho(M)$ is \be\label{Indeterminacy1}\mathbb{L}^{\op{K}}F(X)\doteq F(\tilde{C}_{\tt M}X)\doteq\mathbb{L}_Q^{\op{S}}F(X)=F(Q_{\tt M} X)\cong\mathbb{L}^{\op{S}}F(X)\;,\ee where $\doteq$ is a canonical isomorphism in $\tt Ho(N)$ and $\cong$ a not necessarily canonical isomorphism. Moreover \cite[Proposition 1]{Models}, if $C_{\tt M}X$ is any cofibrant replacement of $X\,,$ there is a canonical $\tt Ho(N)$-isomorphism \be\label{Indeterminacy2} F(\tilde{C}_{\tt M}X)\doteq F(C_{\tt M}X)\;.\ee In view of \eqref{Indeterminacy1} and \eqref{Indeterminacy2} the value of a derived functor at an object is well-defined only up to isomorphism of the target homotopy category. The isomorphism class is independent of the type of derived functor considered, Kan extension or strongly universal, as well as independent of the type of cofibrant F-replacement considered, local or global. Also the choice of another local or another global replacement does not change the isomorphism class. If we compute the value of the derived functor using a local cofibrant replacement that is not necessarily an F-replacement, we get again the same class. Finally, the three representatives considered of the value of the derived functor, namely $L_{\tt N}(F(\tilde{C}_{\tt M}X)),$ $L_{\tt N}(F(Q_{\tt M}X))$ and $L_{\tt N}(F(C_{\tt M}X))\,,$ are related by canonical isomorphisms when viewed as objects of $\tt Ho(N)$ and by zigzags of weak equivalences when viewed as objects of $\tt N\,$. All of this, of course, is consistent with the idea underlying homotopy theory that instead of requiring two objects to be equal, we should simply ask that they are related by a weak equivalence.

\begin{rem}
\emph{Let us stress that the symbols $\sim$ (resp., $\approx$, $\simeq$, $\cong$, $\doteq$) denote in this text a weak equivalence (resp., a zigzag of weak equivalences, a homotopy, an isomorphism, a canonical isomorphism).}
\end{rem}

\section{Derived tensor product in modules over a DG $\cD$-algebra}\label{DTPF}

It is well-known that if we fix one argument in a bifunctor, i.e., in a functor $\zP:\tt C\times D\to E$ out of a product category, we get a functor in the other argument. More precisely, for every $c\in\tt C$ and $d\in\tt D\,,$ the partial functors $L_c:= \zP(c,-):\tt D\to E$ and $R_d:=\zP(-,d):\tt C\to E$ satisfy \be\label{BiFun1}L_c(d)=\zP(c,d)=R_d(c)\;,\ee and for every morphisms $f:c'\to c''$ and $g:d'\to d''$, we have \be\label{BiFun2}L_{c''}(g)\circ R_{d'}(f)=\zP(\id_{c''},g)\circ \zP(f,\id_{d'})=\zP(f,g)=\zP(f,\id_{d''})\circ\zP(\id_{c'},g)=R_{d''}(f)\circ L_{c'}(g)\;.\ee Conversely, if for every object $(c,d)$ of $\tt C\times D$ there are functors $L_c$ and $R_d$ that satisfy \eqref{BiFun1} and if for every morphism $(f,g)$ of $\tt C\times D$ the condition \eqref{BiFun2} is satisfied, then there exists a bifunctor $\zP$ whose partial functors are the $L_c$ and $R_d$ \cite{MacL}. We use this connection between a bifunctor and its partial functors below.\medskip

From here on we work in the {\small HA} context $(\tt{DG\cD M,DG\cD M, DG\cD A})$ (see Section \ref{Intro}) and more specifically in the symmetric monoidal model category ${\tt Mod}(\cA):={\tt Mod}_{\tt DG\cD M}(\cA)$ of modules in $\tt DG\cD M$ over an arbitrary $\cA\in\tt DG\cD A$ \cite{HAC}. We denote the symmetric monoidal structure by $({\tt Mod}(\cA),\0_\cA,\cA)\,$.\medskip

In view of the {\small HA} properties \cite{HAC}, the universal arrow $f\square g$ in the next pushout-diagram is a cofibration if $f:M\to N$ and $g:P\to S$ are and it is a trivial cofibration if one of $f$ and $g$ is a cofibration and the other a trivial cofibration:
$$
\xymatrix@R=1cm@C=1cm{
 M\otimes_\cA P \ar[r]^{f\otimes  \op{id}}      \ar[d]_{\id\otimes g }         & N\otimes_\cA P \ar[d] \ar@/^1pc/[rdd]^{\op{id}\otimes g}   \\
 M\otimes_\cA S \ar[r] \ar@/_1pc/[rrd]_{f\otimes \op{id}} & (N\otimes_\cA P)\, \coprod_{\,M\otimes_\cA P}\,(M\otimes_\cA S)  \ar@{-->}[rd]^{f\square g}             \\
                                      &                                   & N\otimes_\cA S}
$$
If $U\in{\tt Mod}(\cA)$ is cofibrant, the unique arrow $f:\{0\}\to U$ from the initial $\cA$-module to $U$ is a cofibration, so the universal arrow $$f\square g=\id_U\0_\cA\, g:U\otimes_\cA P\to U\0_\cA S$$ is a cofibration (resp., a trivial cofibration), if $g:P\to S$ is a cofibration (resp., a trivial cofibration). Hence, the functor $U\0_\cA-$ respects cofibrations and trivial cofibrations and is therefore a left Quillen endofunctor of the model category ${\tt Mod}(\cA)\,.$ Given what we said above, the functor \be\label{WE1}U\0_\cA -:{\tt Mod}(\cA)\to{\tt Mod}(\cA)\ee sends weak equivalences between cofibrant objects to weak equivalences and of course so does the functor \be\label{WE2}-\0_\cA\, U:{\tt Mod}(\cA)\to{\tt Mod}(\cA)\;.\ee Now, if $(h,k): (U,U')\stackrel{\sim}{\to} (V,V')$ is a weak equivalence between cofibrant objects of the product model category ${\tt Mod}(\cA)\times{\tt Mod}(\cA),$ i.e., if $(h,k)$ is a pair of weak equivalences between cofibrant objects of ${\tt Mod}(\cA)\,,$ then in view of \eqref{BiFun2}, \eqref{WE1} and \eqref{WE2}, the ${\tt Mod}(\cA)$-morphism
$$
h\0_\cA k = (V\0_\cA k)\circ(h\0_\cA U')
$$
is a weak equivalence. Hence the functor $$-\0_\cA -:{\tt Mod}(\cA)\times{\tt Mod}(\cA)\to{\tt Mod}(\cA)$$ meets the requirements of Theorem \ref{Fundamental0} and all reasonable definitions of its left derived functor lead to isomorphic results.

\begin{rem}\label{Fix}
\emph{Let us} fix a functorial Cof - TrivFib factorization \emph{system of $\,{\tt Mod}(\cA)\,,$ so that the cofibrant F-replacement functor $Q$ and the natural weak equivalence $q:Q\stackrel{\sim}\Rightarrow \id$ are fixed as well. In the following, we write} $$-\0_\cA^{\mathbb{L}}-:=\mathbb{L}^{\op{S}}_{Q\times Q}(-\0_\cA -)\;.$$
\end{rem}

However, the left Quillen endofunctor $U\0_\cA -$ ($U$: cofibrant) not only sends weak equivalences between cofibrant objects to weak equivalences, but it preserves {\it all} weak equivalences (\cite[Proposition 3.4.1.]{HAC}). Even better:

\begin{prop}\label{TPPWEQC}
Let $M,P,S\in\tt Mod(\cA)$ and let $g:P\to S$ be a weak equivalence.
\begin{enumerate}
\item If $M$ is cofibrant, the map $\id_M\0\, g:M\0_\cA P\to M\0_\cA S$ is a weak equivalence.
\item If $P$ and $S$ are cofibrant, the map $\id_M\0\,g:M\0_\cA P\to M\0_\cA S$ is a weak equivalence.
\end{enumerate}
\end{prop}

\begin{proof}
We must still prove Item 2. In view of Item 1 and Equation \eqref{BiFun2}, we have the commutative diagram
\begin{center}
\begin{tikzcd}
QM\otimes_\cA P\arrow[r,"\sim"]\arrow[d,"\sim"]& QM\otimes_\cA S\arrow[r,"\sim"]& M\otimes_\cA S\\
M\otimes_\cA P\arrow[rru]       &   &
\end{tikzcd}
\end{center}
From the 2-out-of-3 property it follows that $M\otimes_\cA P\to M\otimes_\cA S$ is a weak equivalence as well.
\end{proof}

The previous proposition implies that in our model categorical context it suffices, just like in classical homological algebra, to resolve one of the two arguments:

\begin{theo}\label{DerTenCRS}
The composites $L \star\0_\cA\star(q\times\id_Q)$ and $L \star\0_\cA\star(\id_Q\times q)\,,$ where $\star$ denotes whiskering and $\id_Q$ is the natural automorphism of $\,Q\,,$ are canonical natural isomorphisms
\be\label{DerF4}{\tt Ho}(L\circ(\id_{{\tt Mod}(\cA)}-\0_\cA\,Q-))\doteq -\0_\cA^{\mathbb{L}}-\,:={\tt Ho}(L\circ(Q-\0_\cA\,Q-))\doteq{\tt Ho}(L\circ(Q-\0_\cA\,\id_{{\tt Mod}(\cA)}-))\;.\ee
Hence, for every $M,N\in{\tt Mod}(\cA)\,,$ we have \be\label{DerF5}M\0_\cA^{\mathbb{L}}N\,:=QM\0_\cA\,QN\approx QM\0_\cA\,N\approx M\0_\cA\,QN\;\ee in ${\tt Mod}(\cA)\,.$ Further: \be\label{DerF6}M\0_\cA^{\mathbb{L}}N\,:=QM\0_\cA\,QN\approx \tilde{C}M\0_\cA\tilde{C}N\approx CM\0_\cA CN\approx QM\0_\cA CN\;,\ee where $\tilde{C}$ (resp., $C$) denotes any local cofibrant F-replacement (resp., any local cofibrant replacement).
\end{theo}

\begin{proof}
Of course, it suffices to prove one of the two statements of \eqref{DerF4}, for example the second. If $f:M\stackrel{\sim}{\to}N$ and $g:P\stackrel{\sim}{\to}S$ are two weak equivalences in ${\tt Mod}(\cA)\,,$ it follows from $$ Qf \0_\cA g = (Qf \0_\cA S)\circ(QM \0_\cA g)$$ and Proposition \ref{TPPWEQC} that the functor $L\circ(Q-\0_\cA\,\id_{{\tt Mod}(\cA)}-)$ sends every weak equivalence of ${\tt Mod}(\cA)^{\times 2}$ to an isomorphism of ${\tt Ho}({\tt Mod}(\cA))\,,$ so that its factorization ${\tt Ho}(L\circ(Q-\0_\cA\,\id_{{\tt Mod}(\cA)}-))$ is well-defined. Due to \cite[Lemma 1]{CompTheo} we only need to build a canonical natural isomorphism $$L\circ(Q-\0_\cA\,Q-)\doteq L\circ(Q-\0_\cA\,\id_{{\tt Mod}(\cA)}-)\;.$$ Since $\id_Q:Q\Rightarrow Q$ is a natural automorphism and $q:Q\Rightarrow\id_{{\tt Mod}(\cA)}$ is a natural weak equivalence, the horizontal composite $\zy:=L\star \0_\cA\star(\id_Q\times q)$ is a canonical natural transformation $$\zy:L\circ(Q-\0_\cA\,Q-)\Rightarrow L\circ(Q-\0_\cA\,\id_{{\tt Mod}(\cA)}-)$$ whose components $\zy_{M,N}=L(QM\0_\cA q_N))$ are isomorphisms. Equation \eqref{DerF5} is a direct consequence of \cite[Proposition 2]{Models} and Equation \eqref{DerF6} follows from Theorem \ref{Fundamental0} and Equation \eqref{Indeterminacy2}.
\end{proof}

\section{Graded $\op{Tor}$ functors in modules over a graded $\cD$-algebra}\label{GraTor}

In contrast to the preceding section, we will mainly work here in the category ${\tt Mod}(\mathfrak{A}):={\tt Mod}_{\tt G\cD M}(\mathfrak{A})$ of modules in the category $\tt G\cD M$ over $\mathfrak{A}\in\tt G\cD A\,$. There are no noteworthy differences to the similar category ${\tt Mod}(\cA):={\tt Mod}_{\tt DG\cD M}(\cA)$ ($\cA\in\tt DG\cD A$) that we considered above. In particular the category ${\tt Mod}(\mathfrak{A})$ is symmetric monoidal with tensor product $-\0_{\mathfrak{A}}-$ which is defined similarly to $-\0_\cA-\,.$

\begin{rem}\label{Not} \emph{It is well-known \cite{BD04} that the category ${\tt Mod}_{\tt \cD M}(\text{\sf{A}})$ ($\text{\sf{\cA}}\in\tt \cD A$) coincides with the ca\-te\-go\-ry $\text{\sf{A}}[\cD]{\tt M}$ of modules over the ring $\text{\sf{A}}[\cD]:=\sf{A}\0_\cO\cD$ (tensor product over functions) of linear differential operators with coefficients in $\text{\sf{A}}\,,$ whose multiplication is defined for instance in \cite{PP}. Moreover, the category ${\tt Mod}_{\tt DG\cD M}(\text{\sf{A}})$ coincides with the category $\tt DG\text{\sf{A}}[\cD]M$ of (non-negatively) graded chain complexes in $\text{\sf{A}}[\cD]{\tt M}$. However, the fact that $\text{\sf{A}}\in\cD\tt A$ is not graded is crucial in the proof of these results. Therefore we will not use notations like $\tt G\mathfrak{A}[\cD]M$ and $\tt DG\cA[\cD]M$ for the categories $\tt Mod(\mathfrak{A})$ and $\tt Mod(\cA)\,,$ respectively.}
\end{rem}

The category ${\tt Mod}_{\tt\cD M}(\text{\sf A})$ is a category of modules over a ring and is therefore abelian and has enough projectives (at least if we assume the axiom of choice, what we do systematically). Further:

\begin{prop}The categories $\tt Mod(\mathfrak{A})$ and $\tt Mod(\cA)$ are abelian and $\tt Mod(\mathfrak{A})$ has enough projectives.\end{prop}

\begin{proof} The categories $\tt G\cD M$ and $\tt DG\cD M$ are abelian symmetric monoidal categories (see The Stacks project, Section 12.16.2 and see \cite{KTRCR}). It therefore follows from \cite{Ard} that the categories ${\tt Mod}_{\tt G\cD M}(\mathfrak{A})$ and ${\tt Mod}_{\tt DG\cD M}(\cA)$ are also abelian if the functors $$\mathfrak{A}\0_\cO-:\tt G\cD M\to G\cD M\quad\text{and}\quad\cA\0_\cO-:\tt DG\cD M\to DG\cD M$$ are additive and preserve cokernels. As the $\cO$-module isomorphisms $$M\0_\cO\{0\}\cong \{0\}\quad\text{and}\quad M\0_\cO(N\oplus P)\cong M\0_\cO N\oplus M\0_\cO P$$ respect the gradings and the actions of vector fields, as well as the differentials if they are present, both left tensoring functors are additive. As for cokernels, notice that both tensoring functors are right exact in $\cO$-modules, hence commute with all finite colimits, in particular with cokernels in $\cO$-modules. However, cokernels in (differential) graded $\cD$-modules are computed degree-wise in $\cD$-modules and cokernels in $\cD$-modules coincide with cokernels in $\cO$-modules and are just the quotients by the set-theoretical images. From here it follows that the functors $\mathfrak{A}\0_\cO-$ and $\cA\0_\cO-$ commute with cokernels in $\tt G\cD M$ and $\tt DG\cD M$, respectively, and that the categories $\tt Mod(\mathfrak{A})$ and $\tt Mod(\cA)$ are abelian.\medskip

Before we go any further, we recall some general definitions.

\begin{defi}\begin{enumerate}\item An object $C$ of a category $\tt C$ is \emph{projective} if it has the {\small LLP} with respect to epimorphisms, or, equivalently, if its covariant Hom-functor $\h(C,-)$ respects epimorphisms. \item A category $\tt C$ has \emph{enough projectives} if for any $C\in\tt C$ there is an epimorphism $P\to C$ from a projective object $P\,.$ \item A \emph{projective resolution} of an object $A$ of an abelian category $\tt A$ is a chain complex $P_\bullet\in{\tt Ch}_+({\tt A})$ made of projective objects $P_i$ together with a quasi-isomorphism $P_\bullet\to A$, or, equivalently, it is an exact sequence $$\cdots \to P_1\to P_0\to A\to 0$$ with projective nods $P_i\,.$\end{enumerate}\end{defi}

\begin{rem}
\emph{The projective resolutions of $A\in{\tt A}$ are exactly the cofibrant replacements of $A\in{\tt Ch}_+({\tt A})$ in the projective model structure of ${\tt Ch}_+({\tt A})\,.$ If $\tt A$ has enough projectives any of its objects has a projective resolution.}
\end{rem}

We still have to prove that $\tt Mod(\mathfrak{A})$ has enough projectives. It is clear that the projective objects in $\tt G\mathcal{D}M$ are exactly the direct sums over $\N$ of projective objects in $\tt \cD M$ (projective objects in $\tt DG\cD M$ are more complicated) and that the category $\tt\cD M$ of modules over $\cD$ has enough projectives. If $M=\oplus_i M_i\in\tt G\cD M\,,$ there exists for each $M_i\in\tt\cD M$ a surjective $\cD$-linear map $P_i\to M_i$ out of a projective $P_i\in\tt \cD M\,,$ hence, there exists a surjective degree-respecting $\cD$-linear map $h:P\to M$ out of the projective $P:=\oplus_i P_i\in\tt G\cD M\,,$ so that also $\tt G\cD M$ has enough projectives. To see that the same holds for ${\tt Mod}(\mathfrak{A})\,$, recall the free-forgetful adjunction $$\mathfrak{A}\0_\cO-:{\tt G\cD M}\rightleftarrows{\tt Mod}(\mathfrak{A}):\op{For}\;,$$ i.e., the functorial bijections $$\Phi_{N M}:\h_{{\tt Mod}(\mathfrak{A})}(\mathfrak{A}\0_\cO N,M)\cong\h_{\tt G\cD M}(N,\op{For}M)\;,$$ which are rooted in the observation that $\mathfrak{A}$-linear maps $$k(a\0 n)=a\triangleleft\, k(1_{\mathfrak{A}}\0 n)$$ are fully defined by their values on the elements $1_{\mathfrak{A}}\0 n\cong n\in N\,.$ From here it follows that for $M\in{\tt Mod}(\mathfrak{A})$ the above surjective $\tt G\cD M$-map $h:P\to M$ is the image by $\Phi$ of a unique ${\tt Mod}(\mathfrak{A})$-map $k:\mathfrak{A}\0_\cO P\to M\,.$ As $k$ is obviously surjective, we can conclude that ${\tt Mod}(\mathfrak{A})$ has enough projectives as soon as we have shown that $\mathfrak{A}\otimes_\cO P\in {\tt Mod}(\mathfrak{A})$ is projective since $P\in\tt G\cD M$ is. However, this is obvious. Indeed, if we choose in ${\tt Mod}(\mathfrak{A})$ a surjective morphism $f:R\to S$ and a morphism $g_{\mathfrak{A}}:\mathfrak{A}\0_\cO P\to S\,,$ its image $g:=\Phi(g_{\mathfrak{A}}):P\to S$ lifts to $R\,,$ i.e., there is a $\tt G\cD M$-morphism $\ell:P\to R$ such that $f\circ \ell= g\,.$ The image $\ell_{\mathfrak{A}}:=\Phi^{-1}(\ell):\mathfrak{A}\0_\cO P\to R$ then lifts $g_{\mathfrak{A}}\,,$ i.e., $f\circ\ell_{\mathfrak{A}}=g_{\mathfrak{A}}\,.$
\end{proof}

For $M\in{\tt Mod}(\mathfrak{A})\,,$ the tensor product $M\0_\mathfrak{A}-:{\tt Mod}(\mathfrak{A})\to{\tt Mod}(\mathfrak{A})$ is a covariant functor between abelian categories, whose source has enough projectives. We can therefore consider its classical left derived functors $L_p(M\0_\mathfrak{A}-):{\tt Mod}(\mathfrak{A})\to {\tt Mod}(\mathfrak{A})$ ($p\in\N$), if the tensor product is right exact.

\begin{rem}
\emph{In Section \ref{DTPF}, we defined the left derived functor $M\!\0_\cA^{\mathbb{L}}-$ (in the model categorical sense) of the tensor product in ${\tt Mod}(\cA)\,,$ for $\cA\in\tt DG\cD A\,.$ In the present section \ref{GraTor}, we are considering the left derived functors $L_p(M\0_\mathfrak{A}-)$ (in the classical homological algebraic sense) of the tensor product in ${\tt Mod}(\mathfrak{A})\,,$ for $\mathfrak{A}\in\tt G\cD A\,.$ We will define the graded Tor functor using the functors $L_p(M\0_\mathfrak{A}-)\,,$ and the Tor spectral sequence will connect the graded Tor functor with the derived functor $M\!\0_\cA^{\mathbb{L}}-\,.$}
\end{rem}

\begin{rem}\label{sest}
\emph{Let us also stress the so far implicit fact that in this text we consider left module structures, left $\cD$-modules and left modules over the graded-commutative algebras $\cA$ and $\mathfrak{A}$ -- although a left action implements a right action via $m\triangleright a=(-1)^{am}a\triangleleft m$ (the numbers $a$ and $m$ in the exponent of $-1$ are the degrees of the vectors $a$ and $m$) and vice versa.}
\end{rem}

\begin{prop}
The tensor product functor $$M\0_{\mathfrak{A}}-:\tt Mod(\mathfrak{A})\to Mod(\mathfrak{A})$$ is right exact for every $M\in\tt Mod(\mathfrak{A})\,.$ The same is true for $-\0_\mathfrak{A}M\,.$
\end{prop}

\begin{proof}
Notice first that the cokernels, kernels and images of the morphisms $\ell:M\to N$ of ${\tt Mod}(\mathfrak{A})$ are those of the underlying morphisms of $\tt G\cD M$ (i.e., the direct sums of the cokernels, kernels and images of the component morphisms of $\cD$-modules). Indeed, the set-theoretical image $\op{im}\ell\in\tt G\cD M$ is an $\mathfrak{A}$-submodule of $N,$ so that the cokernel $\op{coker}\ell =N/\op{im}\ell\in\tt G\cD M$ is also an $\mathfrak{A}$-module. It is easily seen that $\op{coker}\ell\in{\tt Mod}(\mathfrak{A})$ is the cokernel of $\ell$ in ${\tt Mod}(\mathfrak{A})\,.$ Similarly, kernels in ${\tt Mod}(\mathfrak{A})$ are given by the corresponding set-theoretical kernels in $\tt G\cD M$ equipped with their induced $\mathfrak{A}$-module structure. The same statement is therefore also valid for images.\medskip

We are now ready to prove that $M\0_\mathfrak{A}-$ is a right exact endofunctor of $\tt Mod(\mathfrak{A})\,.$ Let therefore \be\label{Seq1}P\stackrel{f}{\to} R\stackrel{g}{\to} S\to 0\ee be an exact sequence in the abelian category ${\tt Mod}(\mathfrak{A})\,,$ i.e., a sequence in ${\tt Mod}(\mathfrak{A})$ such that at each spot the image of the incoming map coincides with the kernel of the outgoing one, i.e., in view of what we just said, a sequence such that at each spot the set-theoretical image coincides with the set-theoretical kernel.\medskip

Observe now that $\mathfrak{A}\in\tt G\cD A$ is a graded-commutative ring $\tilde{\mathfrak{A}}$ (all the rings considered in this text are implicitly assumed to be unital) when equipped with its additive group structure and its multiplication. Any object $P\in \tt Mod(\mathfrak{A})$ is in particular an additive group that comes equipped with a degree zero, $\cD$-linear $\zn:\mathfrak{A}\0_\cO P\to P$ that satisfies the usual associativity and unitality requirements. Since $\zn$ is $\cO$-linear on $\0_\cO$, it is $\cO$-bilinear on $\times$, so in particular biadditive. Hence $\zn$ is an action of the ring $\tilde{\mathfrak{A}}$ on $P\,,$ which is therefore an $\tilde{\mathfrak{A}}$-module $\tilde{P}\,:$ $\tilde P\in\tt Mod(\tilde{\mathfrak{A}})\,.$ Finally a $\tt Mod(\mathfrak{A})$-morphism $f:P$ $\to R$ is a degree zero, $\cD$- and $\mathfrak{A}$-linear map, so it is in particular a ${\tt Mod}(\tilde{\mathfrak{A}})$-morphism $\tilde{f}:\tilde{P}\to\tilde{R}\,.$ From what we just said follows that the exact $\tt Mod(\mathfrak{A})$-sequence \eqref{Seq1} can be interpreted as an exact $\tt Mod(\tilde{\mathfrak{A}})$-sequence \be\label{Seq2}\tilde{P}\stackrel{\tilde{f}}{\to} \tilde{R}\stackrel{\tilde{g}}{\to}\tilde{S}\to 0\;.\ee Actually, the $\tt Mod(\mathfrak{A})$-sequence \eqref{Seq1} is exact if and only if the corresponding $\tt Mod(\tilde{\mathfrak{A}})$-sequence \eqref{Seq2} is.\medskip

We have to show that $M\0_\mathfrak{A}-$ respects exactness of \eqref{Seq1}, i.e., that the $\tt Mod(\mathfrak{A})$-se\-quence \be\label{Seq3}M\0_\mathfrak{A}P\stackrel{M\0 f}{\longrightarrow} M\0_\mathfrak{A}R\stackrel{M\0 g}{\longrightarrow} M\0_\mathfrak{A}S\longrightarrow 0\ee is exact, or, equivalently, that the corresponding $\tt Mod(\tilde{\mathfrak{A}})$-sequence \be\label{Seq3}(M\0_\mathfrak{A}P)^{\widetilde{}}\,\stackrel{(M\0 f)^{\tilde{}}}{\longrightarrow} \,(M\0_\mathfrak{A}R)^{\widetilde{}}\,\stackrel{(M\0 g)^{\tilde{}}}{\longrightarrow}\,(M\0_\mathfrak{A}S)^{\widetilde{}}\,\longrightarrow\, 0\ee is, i.e., that at each of its spots the set-theoretical image coincides with the set-theoretical kernel. However, the tensor product $-\0_{\tilde{\mathfrak{A}}}-$ in ${\tt Mod}(\tilde{\mathfrak{A}})$ over the graded-commutative ring $\tilde{\mathfrak{A}}$ being the standard tensor product over the non-commutative ring $\tilde{\mathfrak{A}}\,,$ the functor $\tilde{M}\0_{\tilde{\mathfrak{A}}}-$ is right exact, so that the se\-quence \be\label{Seq4}\tilde{M}\0_{\tilde{\mathfrak{A}}}\tilde{P}\stackrel{\tilde{M}\0\tilde{f}}{\longrightarrow} \tilde{M}\0_{\tilde{\mathfrak{A}}}\tilde{R}\stackrel{\tilde{M}\0\tilde{g}}{\longrightarrow}\tilde{M}\0_{\tilde{\mathfrak{A}}}\tilde{S}\to 0\;\ee is exact in $\tt Mod(\tilde{\mathfrak{A}})\,,$ i.e., at each spot the set-theoretical image coincides with the set-theoretical kernel. It is now sufficient to explain why Sequence \eqref{Seq3} and Sequence \eqref{Seq4} are made of the same sets and the same set-theoretical maps. The maps have obviously both the same values as $M\0 f$ and $M\0 g$, respectively. It is therefore enough to prove that $(M\0_\mathfrak{A}P)^{\widetilde{}}=\tilde{M}\0_{\tilde{\mathfrak{A}}}\tilde P\,.$ By definition the latter set is made of the finite sums of decomposed tensors $m\0 p\,,$ where $\0$ is weakly $\mathfrak{A}$-bilinear, i.e., is biadditive and such that \be\label{WABL}(a\triangleleft m)\0 p=(-1)^{am}m\0 (a\triangleleft p)\;,\ee where we denoted the vectors by the same symbols as their degree. The first tensor product $M\0_\mathfrak{A}P$ is defined as the cokernel in $\tt G\cD M$ of the map $$\za:M\0_\cO\mathfrak{A}\0_\cO P\ni (m,a,p)\mapsto (a\triangleleft m)\0 p - (-1)^{am}m\0(a\triangleleft p)\in M\0_\cO P\;.$$ As said above $\op{coker}\za\in\tt G\cD M$ is the quotient $(M\0_\cO P)/\op{im}\za\,,$ so it is given by the finite sums of decomposed tensors $m\0 p\,,$ where $\0$ is $\cO$-bilinear and weakly $\mathfrak{A}$-bilinear. This graded $\cD$-module inherits a compatible $\mathfrak{A}$-action $$a\triangleleft(m\0 p)=(a\triangleleft m)\0 p=(-1)^{am}m\0(a\triangleleft p)\;.$$ When passing to the $\tilde{\mathfrak{A}}$-module $(M\0_\mathfrak{A}P)^{\widetilde{}}\,,$ we forget in particular the $\cO$-action, so this set is made of the finite sums of decomposed tensors $m\0 p\,,$ where $\0$ is weakly $\mathfrak{A}$-bilinear: $(M\0_\mathfrak{A}P)^{\widetilde{}}=\tilde{M}\0_{\tilde{\mathfrak{A}}}\tilde{P}\,.$\end{proof}

\begin{defi} Let $\mathfrak{A}\in\tt G\cD A$ and $M\in\tt Mod(\mathfrak{A})\,.$ We refer to the $p$-th left derived functor of the covariant right exact endofunctor $M\0_\mathfrak{A}-$ as the $p$-th \emph{graded Tor endofunctor} $\op{Tor}^{\mathfrak{A}}_p(M,-)$ of $\,\tt Mod(\mathfrak{A})$. For any $N\in\tt Mod(\mathfrak{A})\,,$ we have $\op{Tor}^\mathfrak{A}_p(M,N)\in{\tt Mod}(\mathfrak{A})\subset \tt G\cD M$ and we denote the $q$-th homogeneous component by $\op{Tor}^\mathfrak{A}_p(M,N)_q\in\tt\cD M\,$.
\end{defi}

Let us recall that in order to compute $\op{Tor}^\mathfrak{A}_p(M,N)\,,$ we choose a projective resolution $P_\bullet\to N$ of $N$ in the abelian category $\tt Mod(\mathfrak{A})\,,$ i.e., we choose a cofibrant replacement $\mathcal{Q}N$ of $N$ in the projective model structure of ${\tt Ch}_+(\tt Mod(\mathfrak{A}))$ (in particular the cofibrant replacement $\mathcal{Q}N$ that is provided by a functorial `cofibration -- trivial fibration' factorization), and we compute the $p$-th homology $\mathfrak{A}$-module of $M\0_\mathfrak{A}\mathcal{Q}N\in{\tt Ch}_+(\tt Mod(\mathfrak{A}))\,.$ In other words, we have $$\op{Tor}^\mathfrak{A}_p(M,-):{\tt Mod}(\mathfrak{A})\stackrel{i}{\longrightarrow}{\tt Ch}_+({\tt Mod}(\mathfrak{A}))\stackrel{M\0\mathcal{Q}-}{\longrightarrow}{\tt Ch}_+({\tt Mod}(\mathfrak{A}))\stackrel{H_p}{\longrightarrow}\tt Mod(\mathfrak{A})\;.$$ Here the first arrow $i$ is the functor that sends a module (resp., a module morphism) to the corresponding chain complex concentrated in degree zero (resp., the corresponding chain map concentrated in degree zero), the second arrow $M\0\mathcal{Q}-$ is the composite of a cofibrant replacement functor in the projective model structure of chain complexes and the tensor product functor, whereas the last arrow $H_p$ is the $p$-th homology functor. Finally, as in the classical situation, if we resolve $M$ instead of $N$, we get the same result.\medskip

In view of Theorem \ref{DerTenCRS}, we have in $\tt Mod(\cA)$ ($\cA\in\tt DG\cD A$) that \be\label{L}M\0^{\mathbb{L}}_\cA N:=QM\0_\cA QN\approx QM\0_\cA N\approx M\0_\cA QN\in{\tt Mod}(\cA)\subset{\tt DG\cD M}={\tt Ch}_+(\cD\tt M)\;\ee ($Q$: cofibrant replacement functor in $\tt Mod(\cA)$), whereas in $\tt Mod(\mathfrak{A})$ ($\mathfrak{A}\in\tt G\cD A$) we have that \be\label{T}\op{Tor}^\mathfrak{A}_p(M,N)=H_p(M\0_\mathfrak{A}\mathcal{Q}N)\in{\tt Mod}(\mathfrak{A})\subset{\tt G\cD M}\;\ee ($\mathcal{Q}$: cofibrant replacement functor in ${\tt Ch}_+(\tt Mod(\mathfrak{A}))$). Notice further that, since $M\in\tt Mod(\cA)\subset {\tt DG\cD M}={\tt Ch}_+(\tt\cD M)$ is not a chain complex of $\cA[\cD]$-modules (see Remark \ref{Not}), the homology module $H(M)$ is a priori just a graded $\cD$-module. However, the homology module $H(\cA)\in\tt G\cD M$ inherits the obvious multiplication, which makes it a graded $\cD$-algebra $\mathfrak{A}$. It is also easily seen that the homology module $H(M)\in\tt G\cD M$ can be endowed with the canonical $H(\cA)$-action, which makes it a module $$H(M)\in{\tt Mod}(H(\cA))={\tt Mod}_{\tt G\cD M}(H(\cA))\;;$$ in fact, the functor $H:\tt DG\cD M\to G\cD M$ restricts to a functor \be\label{HomFunModA}H:{\tt Mod}(\cA)\to{\tt Mod}(H(\cA))\;.\ee A possible relationship between $\0_\cA^{\mathbb{L}}$ in \eqref{L} and $\op{Tor}^{\mathfrak{A}}$ in \eqref{T} should therefore involve $$\op{Tor}^{H(\cA)}_p(H(M),H(N))_q\in{\tt \cD M}\quad\text{and}\quad H_{p+q}(M\0_\cA^{\mathbb{L}}N)\in\cD\tt M\;,$$ where $\cA\in\tt DG\cD A\,,$ $M,N\in\tt Mod(\cA)$ and $p,q\in\N\,.$ The $\cD$-generalization of Quillen's $\op{Tor}$ spectral sequence specifies this relationship.

\section{The $\mathcal{D}$-generalization of Quillen's $\op{Tor}$ spectral sequence}\label{TSSDGeo}

\begin{theo}\label{DTSS}
For every $\cA\in{\tt DG\cD A}$ and every $M,N\in\tt Mod(\cA):=Mod_{DG\cD M}(\cA)\,,$ there is a first quadrant spectral sequence $E_{pq}^\bullet$ in the abelian category $\tt \cD M\,,$ whose second sheet is $E_{pq}^2=\op{Tor}^{H(\cA)}_p(H(M),H(N))_q$ and which converges to $H_{p+q}(M\0_\cA^{\mathbb{L}}N)\,:$ \be\label{TorSpecSeq}E_{pq}^2=\op{Tor}^{H(\cA)}_p(H(M),H(N))_q\Rightarrow H_{p+q}(M\0_\cA^{\mathbb{L}}N)\;.\ee
\end{theo}

What we call Quillen's $\op{Tor}$ spectral sequence is a similar spectral sequence in the category $\tt Ab$ of abelian groups, for a simplicial ring $\cA$ and a simplicial right (resp., left) $\cA$-module $M$ (resp., $N$) (for simplicial modules, see for instance `The Stacks project, 21.40'; see also \cite{Linear}, where another type of modules over a varying ring appears). Roughly speaking the result connects the derived tensor product of the homology and the homology of the derived tensor product.\medskip

The proof of Theorem \ref{DTSS} relies on two propositions.

\begin{prop}\label{tri} For every $\cA\in\tt DG\mathcal{D}A$, every $M\in \tt DG\mathcal{D}M$ and every $N\in {\tt Mod(\cA)}$, there is a natural isomorphism $$(M\otimes_\cO \cA)\otimes_\cA N\cong M\otimes_\mathcal{O} N\;$$ in $\,{\tt Mod(\cA)},$ where the $\cA$-module structure of $M\0_\cO N$ is canonically induced by the $\cA$-module structure of $N$: $$(\id_M\otimes\,\nu_N)\circ(\op{com}\0\id_N): \cA\0_\cO M\0_\cO N\ni a\0m\0n\mapsto(-1)^{am}m\0(a\triangleleft n)\in M\otimes_\mathcal{O} N\;,$$ where $\zn_N$ is the $\cA$-module structure $\triangleleft$ of $N$ and $\op{com}$ the braiding of $\,\tt DG\cD M\,$. The same result is valid for $\mathfrak{A}\in\tt G\cD A,$ $M\in\tt G\cD M$ and $N\in{\tt Mod}(\mathfrak{A})\,.$
\end{prop}

Similar results exist in other settings. The main difference here is that $\cA\in\tt DG\cD A$ and $N\in\tt DG\cD M\,,$ but the $\tt DG\cD M$-map $\zn_N:\cA\0_\cO N\to N$ is not $\cD$-bilinear as the tensor product is over $\cO\,.$\medskip

For the proof of Proposition \ref{tri}, we refer the reader to Lemma 3.1.1. in \cite{HAC}.\medskip

For the next proposition we need some preparation. Let $\tt GSet$ be the category whose objects are the (non-negatively) graded sets, i.e., the disjoint unions $I=\sqcup_{k\in\N}I_k$ of sets $I_k$ indexed by the non-negative integers $k\in\N\,,$ and whose morphisms are the grading preserving set-theoretical mappings. We systematically and implicitly consider such sets as well-ordered. Namely, we well-order all the sets $I_k$ and then well-order $I$ by declaring that the elements of $I_k$ are smaller than the elements of $I_\ell\,,$ if $k<\ell$. Every $i\in I$ belongs to a unique $I_k\,,$ so that we get the degree assigning map $n:I\ni i\mapsto n_i=k\in\N\,.$ The functors \be\label{AdjGSet}\mathfrak{F}:{\tt GSet}\rightleftarrows{\tt G\cD M}:\op{For}_\cD\ee form a free-forgetful adjunction. More precisely, the free graded $\cD$-module over a graded set $I$ is the free $\cD$-module $$\mathfrak{F}(I)=\oplus_{i\in I}\cD\, i=\oplus_{k\in\N}\oplus_{i_k\in I_k}\cD\, i_k$$ over the set $I,$ but with $i_k$ in degree $n_{i_k}=k\,,$ or, equivalently, is the direct sum over all degrees $k\in\N$ of the free $\cD$-modules over the sets $I_k\,.$ In other words, we set $$\mathfrak{F}(I)=\oplus_{i\in I}\cD\,{\bf 1}_i[n_i]=\oplus_{i\in I} S_i^{n_i}\in{\tt G\cD M}\;,$$ where ${\bf 1}[n_i]$ is the generator of the $n_i$-sphere $S^{n_i}$ (the $n_i$-sphere is concentrated on $\cD$ in degree $n_i\in\N$). When composing, for $\mathfrak{A}\in\tt G\cD A\,,$ the adjunction \eqref{AdjGSet} with the adjunction $$\mathfrak{A}\0_\cO-:{\tt G\cD M}\rightleftarrows{\tt Mod}(\mathfrak{A}):\op{For}_{\mathfrak{A}}\;,$$ we see that the free $\mathfrak{A}$-module in $\tt G\cD M$ over a graded set $I$ with degree assigning map $n$ is $$\mathfrak{A}\0_\cO(\oplus_{i\in I}S_i^{n_i})=\mathfrak{A}\triangleleft(1_{\mathfrak{A}}\0_\cO(\oplus_{i\in I} S_i^{n_i}))\in{\tt Mod}(\mathfrak{A})\;.$$ Hence:

\begin{lem}\label{DefMorFreeGra}
If $P\in\tt Mod(\mathfrak{A})\,,$ a $\tt Mod(\mathfrak{A})$-morphism \be\label{Basis1}q:\mathfrak{A}\0_\cO(\oplus_{i\in I}S_i^{n_i})\to P\ee is uniquely defined by its values $$q(1_{\mathfrak{A}}\0 {\bf 1}_i[n_i])\in P_{n_i}\quad(i\in I)\;.$$
\end{lem}

A similar result holds in the differential graded setting. More precisely, from item (i) of Lemma 2.3.1. in \cite{HAC} it follows that, if $\cA\in\tt DG\cD A\,,$ the graded $\cD$-module $\cA\0_\cO(\oplus_{i\in I}S_i^{n_i})$ equipped with its natural differential $d_\cA\0\id$ and its natural $\cA$-action is a module in $\tt Mod(\cA)\,.$ From item (ii) of the same lemma we get:

\begin{lem}\label{DefMorFreeDiffGra} If $P\in\tt Mod(\cA)\,,$ a $\tt Mod(\cA)$-morphism \be\label{Basis2}q:\cA\0_\cO(\oplus_{i\in I}S_i^{n_i})\to P\ee is uniquely defined by its values $$q(1_\cA\0{\bf 1}_i[n_i])\in P_{n_i}\cap \ker d_P\quad (i\in I)\;$$ $($provided the source module is equipped with its natural differential $\id_\cA\0\id$$\,)$.
\end{lem}

Indeed, it is then enough to extend $q$ as a $\cD$-linear map to the direct sum and as an $\cA$-linear map to the tensor product.

\begin{prop}\label{lemma2} For $\cA\in \tt DG\mathcal{D}A$ and $P, N\in{{\tt Mod(\cA)}:=\tt Mod}_{\tt DG\cD M}(\cA)\,,$ with $P$ a cofibrant object whose homology $H(P)$ is a free $H(\cA)$-module in $\tt G\cD M$, there exists an isomorphism
$$H(P)\otimes_{H(\cA)}H(N)\cong H(P\otimes_\cA N)\;$$
in the category ${\tt Mod}(H(\cA)):={\tt Mod_{G\cD M}}(H(\cA))\,$.
\end{prop}

\begin{proof} It is natural to interpret the graded $\cD$-module $$S:=\oplus_{i\in I} S_i^{n_i}$$ as a differential graded $\cD$-module with zero differential (since each term is such a differential graded $\cD$-module -- we have used this fact already above). Hence $S$ is a chain complex of $\cO$-modules that is in each degree $k$ a direct sum $S(k)$ of copies of $\cD\,.$ Since $\cD$ is $\cO$-flat (\cite[Proposition 1.0.11.]{HAC}), this implies that $S(k)$, $d(S(k))=0$ and $H_k(S)=S(k)$ are $\cO$-flat, so that for every $M\in\tt DG\cD M\,,$ we have $H(S\0_\cO M)\cong S\0_\cO H(M)\,,$ in view of Künneth's theorem (the homology does of course not depend on whether we interpret a complex as complex of $\cD$-modules or as complex of $\cO$-modules). In particular:

\begin{lem}\label{Kun} If $\cA\in\tt DG\cD A$ and $N\in{\tt Mod}(\cA)\,,$ we have \be\label{KunnethS}S\0_\cO H(N)\cong H(S\0_\cO N)\quad \text{and} \quad H(\cA)\0_\cO S\cong H(\cA\0_\cO S)\;.\ee The first $\,\tt G\cO M$-isomorphism $s\0[n]\mapsto [s\0 n]$ is also a $\,{\tt Mod}(H(\cA))$-isomorphism. A similar statement holds for the second isomorphism.
\end{lem}

Recall now that the cofibrant objects in $\tt Mod(\cA)$ are exactly the retracts of the Sullivan $\cA$-modules, i.e., of the $\cA$-modules $\cA\0_\cO (\oplus_{i\in I}\cD\,g_i[n_i]),$ where $I$ is a well-ordered set (equivalently, an ordinal) and the $g_i[n_i]$ are symbols of degree $n_i\in\N$, and whose differential is lowering in the sense of Definition 2.3.6. of \cite{HAC}.\medskip

By assumption $P$ is a cofibrant $\cA$-module whose homology is $$H(P)=H(\cA)\0_\cO S\cong H(\cA\0_\cO(\oplus_{i\in I}S_i^{n_i}))=H(\cA\0_\cO(\oplus_{i\in I}\cD\,{\bf 1}_i[n_i])\;.$$ As mentioned above, we view $I$ as a well-ordered set; moreover, we have $d({\bf 1}_i[n_i])=0\,,$ so that the differential is obviously lowering. Hence, the $\cA$-module $\cA\0_\cO(\oplus_{i\in I}\cD\,{\bf 1}_i[n_i])$ is a Sullivan module and therefore a cofibrant module.\medskip

In the following, we denote homology classes as usual with $ [-] \,$. Confusion with the degree shift $ [n_i] $ is excluded, since the meaning of $ [-] $ always comes from the context. Since $[1_\cA]\0{\bf 1}_i[n_i]\in H_{n_i}(P)$ for every $i\,,$ we have $[1_\cA]\0{\bf 1}_i[n_i]=[p_i]\,,$ with $p_i\in P_{n_i}\cap\, \ker d_P\,.$ Choosing a representative $p_i$ for each $i\in I$, we get values $$q(1_\cA\0{\bf 1}_i[n_i])=p_i\in P_{n_i}\cap\ker d_P\quad (i\in I)\;,$$ which (see Lemma \eqref{DefMorFreeDiffGra}) define a $\tt Mod(\cA)$-morphism $$q:\cA\0_\cO(\oplus_{i\in I}S_i^{n_i})\to P\;.$$ If $$H(q):H(\cA)\0_\cO(\oplus_{i\in I}S_i^{n_i})\to H(P)$$ denotes the induced ${\tt Mod}(H(\cA))$-morphism in homology (see Equation \eqref{HomFunModA}), we have $$H(q)([1_\cA]\0{\bf 1}_i[n_i]) \cong H(q)[1_\cA\0{\bf 1}_i[n_i]]=[q(1_\cA\0{\bf 1}_i[n_i])]=[p_i]=[1_\cA]\0{\bf 1}_i[n_i]\in H_{n_i}(P)\;$$ ($i\in I$). We know (see Lemma \ref{DefMorFreeGra}) with $\mathfrak{A}=H(\cA)$) that these values define a unique ${\tt Mod}(H(\cA))$-morphism $h:H(\cA)\0_\cO S\to H(P)\,.$ In the situation under consideration this morphism $h$ is identity. Hence $H(q)=\id\,,$ so that $q$ is a weak equivalence between cofibrant $\cA$-modules. Recall that a weak equivalence $f:K\to L$ in $\tt Mod(\cA)$ is a $\tt Mod(\cA)$-morphism whose underlying $\tt DG\cD M$-morphism is a quasi-isomorphism. It is easy to check though that $H(f):H(K)\to H(L)$ is not only an isomorphism in $\tt G\cD M\,,$ but also an isomorphism in ${\tt Mod}(H(\cA))\,.$ From Proposition \ref{TPPWEQC} (with $g = q$) we now get a ${\tt Mod}(H(\cA))$-isomorphism
\begin{equation}\label{nulta}H((\cA\0_\cO S)\0_\cA N)=H((\cA\otimes_\mathcal{O}(\oplus_{i\in I}S_i^{n_i}))\otimes_\cA N)\cong  H(P\otimes_\cA N)\;.\end{equation}
Using Proposition \ref{tri} and Equation \ref{KunnethS}, we get ${\tt Mod}(H(\cA))$-isomorphisms
\begin{equation}\label{prva}H((\cA\0_\cO S)\0_\cA N)\cong H(S\otimes_\mathcal{O} N)\cong S\otimes_\mathcal{O}H(N)\;,\end{equation}
and
\begin{equation}\label{druga}H(P)\0_{H(\cA)} H(N)=(H(\cA)\0_\cO S)\0_{H(\cA)}H(N)\cong S\0_\cO H(N)\;.\end{equation}
Equations \ref{nulta}, \ref{prva} and \ref{druga} finally give the desired result.
\end{proof}

We are now prepared to prove Theorem \ref{DTSS}.\medskip

Let $\cA\in \tt DG\mathcal{D}A$ and $M,N\in {\tt Mod}(\cA)\,.$ Our goal is to prove that there is a first quadrant spectral sequence in the abelian category $\tt \cD M$ such that \be\label{TorSpecSeqInProof}E^2_{pq}=\op{Tor}_{p}^{H(\cA)}(H(M),H(N))_q\Rightarrow H_{p+q}(M\otimes^\mathbb{L}_\cA N)\;,\ee for each $p,q\in\N\,$.\medskip

Recall first that to any first quadrant double complex $(C_{\bullet\bullet}, d^h, d^v)$ in an abelian category, one associates its total complex $(\op{Tot}(C)_\bullet,d)\,,$ whose grading is $\op{Tot}(C)_n=\oplus_{p+q=n} C_{pq}$ and whose differential is $d=d^h+(-1)^pd^v\,.$ The total complex admits two filtrations, the horizontal filtration $F^h_p(\op{Tot}(C))=\oplus_{r\le p} C_{r\bullet}$ and the vertical filtration $F^v_p(\op{Tot}(C))=\oplus_{s\le p} C_{\bullet s}\,.$ The spectral sequence of the resulting horizontally filtered chain complex $(\op{Tot}(C)_\bullet,d,F^h_\bullet)$ (resp., the vertically filtered chain complex $(\op{Tot}(C)_\bullet,d,F^v_\bullet)$) is the horizontal (resp., vertical) spectral sequence ${}^h\!E^\bullet_{pq}$ (resp., ${}^v\!E^\bullet_{pq}$) of the double complex considered. The second sheet of the horizontal (resp., vertical) spectral sequence is $${}^h\!E^2_{pq}=H_p^h(H_q^v(C_{\bullet\bullet}))\quad (\text{resp.,}\;{}^v\!E^2_{pq}=H_p^v(H_q^h(C_{\bullet\bullet})))\;.$$ Moreover, this spectral sequence converges to the homology of the total complex: $${}^h\!E^2_{pq}=H_p^h(H_q^v(C_{\bullet\bullet}))\Rightarrow H_{p+q}(\op{Tot}(C))\quad(\text{resp.,}\;{}^v\!E^2_{pq}=H_p^v(H_q^h(C_{\bullet\bullet}))\Rightarrow H_{p+q}(\op{Tot}(C)))\;.$$ This means that for any $p,q\in\N\,,$ the horizonal spectral sequence ${}^h\!E_{pq}^r$ ($r\in\N$) stabilizes at some $r(p,q)\in\N$ and \be\label{convergence}{}^h\!E^{r(p,q)}_{pq}\cong {}^h G_p(H_{p+q}(\op{Tot}(C)))\;,\ee where the {\small RHS} is the $p$-th term of the graded space that is associated to the filtered space ${}^h\!F_\bullet(H_{p+q}(\op{Tot}(C)))\,,$ whose filtration is induced by the filtration $F_\bullet^h$ of $\op{Tot}(C)_\bullet\,.$ The dual result holds for the vertical spectral sequence and the vertical filtration. If exactly one row or column of the grid ${}^h\!E_{pq}^2$ ($p,q\in\N$) does not vanish, the horizontal spectral sequence collapses at its second sheet and $H_{n}(\op{Tot}(C))$ ($n\in\N$) is the unique non-zero ${}^h\!E_{pq}^2$ such that $p+q=n$ \cite[Definition 5.2.7]{Wei}. Again, the dual result holds for the vertical spectral sequence.\medskip

Also remember (see Section \ref{GraTor}) that in order to compute the {\small LHS} of \eqref{TorSpecSeqInProof}, we compute the $p$-th homology space of the $q$-th term of the tensor product of $H(N)$ with a projective resolution of $H(M)$ in ${\tt Mod}(H(\cA))\,.$ On the other hand, in order to compute the derived tensor product in the {\small RHS} of \eqref{TorSpecSeqInProof} (up to a zigzag of weak equivalences in ${\tt Mod}(\cA)$), we can (see Theorem \ref{DerTenCRS}) tensor over $\cA$ the cofibrant replacement $X:=QM$ of $M$ in $\tt Mod(\cA)$ with {\bf a} cofibrant replacement of $N\,$.\medskip

Since $\tt Mod(\cA)$ is a cofibrantly generated model category \cite[Theorem 2.2.3.]{HAC} (we denote its set of generating cofibrations by $I$), the small object argument gives a functorial factorization of $\tt Mod(\cA)$-morphisms $f:R\to S$ into an $I$-cell $i:R\rightarrowtail \mathfrak{Q}$ and a trivial fibration $p:\mathfrak{Q}\stackrel{\sim}\twoheadrightarrow S\,$. As the specific cofibrations we call $I$-cells are known to be relative Sullivan $\cA$-modules \cite{HAC}, we get in the case of $f:0\to N$ a cofibrant replacement $\mathfrak{Q}$ of $N$ that is a Sullivan $\cA$-module $$\mathfrak{Q}=\cA\0_\cO V:=\cA\0_\cO\oplus_{k\in K}\cD\; g_k$$ (where the set $K$ is well-ordered, the generators $g_k$ have homogeneous degrees $\deg(g_k):=n_k\in\N$ and the differential on $\mathfrak{Q}$ is lowering).\medskip

In the following, we use the cofibrant replacement $X$ of $M$ and the cofibrant replacement $\mathfrak{Q}$ of $N\,$. Further, we will construct a projective resolution of $H(M)\cong H(X)$ in ${\tt Mod}(H(\cA))$ as the sequence induced in homology by a resolution $P_\bullet$ of $X$ in $\tt Mod(\cA)\,.$ The statement \eqref{TorSpecSeqInProof} of Theorem \ref{DTSS} then results from the above machinery for double complexes, applied to the double complex $(P_\bullet\0_\cA\mathfrak{Q})_\bullet\,.$\medskip

We construct the exact sequence $P_\bullet\to X\to 0$ just mentioned, using an iterative process.\medskip

Set $X_0=X$ and fix a family of homogeneous generators $([x_i])_{i\in I_0}$ of $H(X_0)\in{\tt Mod}(H(\cA))\,.$ Homogeneity means that $x_i$ has a homogeneous degree $n_i\in\N\,,$ so that $x_i\in X_{0,n_i}\cap\ker d_{X_0}\,.$ Notice that the number of generators needs not be finite in any way, so that one can in particular choose as generators all the homogeneous elements of $H(X_0)=\oplus_nH_n(X_0)\,.$ Next we define a $\tt Mod(\cA)$-morphism $$q_0:\cA\0_\cO\oplus_{i\in I_0}S_i^{n_i}\to X_0\;.$$ According to Lemma \ref{DefMorFreeDiffGra} such a morphism is uniquely defined by its values $q_0(1_\cA\0 \mathbf{1}_i[n_i])\in X_{0,n_i}\cap \ker d_{X_0}\,,$ in particular by the choice of representatives $x_i$ of the generating homology classes $[x_i]\,.$ Any trivial cofibration - fibration decomposition of $q_0$ leads to a trivial cofibration $\iota_0: \cA\0_\cO S_0\stackrel{\sim}{\rightarrowtail} P_0$ and a fibration $\zp_0:P_0\twoheadrightarrow X_0$ of $\tt Mod(\cA)$ (where $S_0$ is of course a compact notation for $\cA\0_\cO\oplus_{i\in I_0}S_i^{n_i}$). First, the weak equivalence $\iota_0$ of $\tt Mod(\cA)$ induces an isomorphism $H(\iota_0)$ of ${\tt Mod}(H(\cA)).$ Since, due to Lemma \ref{Kun}, we have also a ${\tt Mod}(H(\cA))$-isomorphism $$H(\cA\0_\cO S_0)\simeq H(\cA)\0_\cO S_0\;,$$ we see that $H(P_0)$ is isomorphic as $H(\cA)$-module to the free $H(\cA)$-module $H(\cA)\0_\cO S_0\,.$ Moreover, we already mentioned that modules of the type $\cA\0_\cO \oplus_{i\in I_0}S_i^{n_i}$ are Sullivan and therefore cofibrant $\cA$-modules. Hence, in the composite $$0\rightarrowtail \cA\0_\cO\oplus_{i\in I_0}S_i^{n_i}\stackrel{\iota_0}{\rightarrowtail} P_0$$ is a cofibration and $P_0$ is a cofibrant $\cA$-module. Second, the $\tt Mod(\cA)$-morphism $\zp_0$ is also a fibration in $\tt DG\cD M\,,$ hence it is surjective in positive degrees $n>0\,.$ On the other hand, the induced ${\tt Mod}(H(\cA))$-morphism $H(\zp_0):H(P_0)\to H(X_0)$ is surjective, since any homology class $[x]\in H(X_0)$ can be expressed as finite combination of generators, so that $$[x]=\sum_{i}[a_i]\triangleleft [x_i]=\sum_i[a_i\triangleleft x_i]=H(q_0)[\sum_i a_i\0 \mathbf{1}_i[n_i]]=H(\zp_0)(H(i_0)([\sum_i a_i\0 \mathbf{1}_i[n_i]]))\;.$$ One easily checks that the $\tt Mod(\cA)$-morphism $\zp_0:P_0\to X_0\,,$ which is surjective in positive degrees and surjective in homology, is actually surjective in all degrees.\medskip

In the previous paragraph, we started from a module $X_0\in\tt Mod(\cA)$ and constructed a cofibrant module $P_0\in\tt Mod(\cA)$ whose homology $H(P_0)$ is a free $H(\cA)$-module, as well as a morphism $\zp_0:P_0\to X_0$ that is surjective and induces a surjective morphism in homology. We now set $$X_1:=\ker(\zp_0:P_0\to X_0)\stackrel{k_1}{\hookrightarrow} P_0\,,$$ where $\stackrel{k_{1}}{\hookrightarrow}$ is the canonical inclusion of the $\cA$-submodule $X_{1}$ into the $\cA$-module $P_0\,,$ and we iterate the process of the previous paragraph. This way we get $\tt Mod(\cA)$-morphisms $$\cdots\hookrightarrow P_2\stackrel{\zp_2}{\to}X_2\stackrel{k_2}{\hookrightarrow} P_1\stackrel{\zp_1}{\to}X_1\stackrel{k_1}{\hookrightarrow} P_0\stackrel{\zp_0}{\to}X_0\to 0\;.$$ The $\tt Mod(\cA)$-sequences \be\label{SequenceA}0\to X_{n+1}\stackrel{k_{n+1}}{\hookrightarrow} P_n\stackrel{\zp_n}{\to} X_n\to 0\ee ($n\in\N$) and \be\label{SequenceB}\cdots \to P_2\stackrel{\tilde{\zp}_2}{\to} P_1\stackrel{\tilde\zp_1}{\to} P_0\stackrel{\zp_0}{\to} X\to 0\;\ee are obviously exact.\medskip

As mentioned earlier, we have to check whether the cofibrant resolution $P_\bullet$ of $X=X_0$ in $\tt Mod(\cA)$ (i.e., resolution made of cofibrant $\cA$-modules) induces a projective resolution $H(P_\bullet)$ of $H(X)$ in ${\tt Mod}(H(\cA))\,.$ When applying the homology functor $H:{\tt Mod}(\cA)\to{\tt Mod}(H(\cA))$ to the sequences \eqref{SequenceA} and \eqref{SequenceB}, we get sequences \be\label{SequenceAH}0\longrightarrow H(X_{n+1})\stackrel{H(k_{n+1})}{\longrightarrow} H(P_n)\stackrel{H(\zp_n)}{\longrightarrow} H(X_n)\longrightarrow 0\;\ee and \be\label{SequenceBH}\cdots \longrightarrow H(P_2)\stackrel{H(\tilde\zp_2)}{\longrightarrow} H(P_1)\stackrel{H(\tilde\zp_1)}{\longrightarrow} H(P_0)\stackrel{H(\zp_0)}{\longrightarrow} H(X)\longrightarrow 0\;.\ee As \eqref{SequenceA} is also a short exact sequence in $\tt Ch_+(\cD M)\,,$ it induces an exact $\tt G\cD M$-triangle in homology. Since $H(\zp_n)$ is surjective, the connecting homomorphism vanishes and the morphism $H(k_{n+1})$ is injective (it is actually the canonical inclusion), so that the short ${\tt Mod}(H(\cA))$-sequence \eqref{SequenceAH} is exact. It follows that the long ${\tt Mod}(H(\cA))$-sequence \eqref{SequenceBH} is exact. Indeed, since $$H(\tilde{\zp}_{n+1})=H(k_{n+1})\circ H(\zp_{n+1})\;,$$ the preceding properties imply that $$\ker H(\tilde{\zp}_{n+1})=\ker H(\zp_{n+1})=H(X_{n+2})=\op{im}H(\zp_{n+2})=\op{im}H(\tilde{\zp}_{n+2})\;.$$ Hence $H(P_\bullet)$ is a free (hence projective) resolution of $H(X)$ in ${\tt Mod}(H(\cA))\,.$\medskip

We continue to follow the procedure described above and consider the first quadrant double complex $(P_p\0_\cA\mathfrak{Q})_q\in\tt\cD M\,$ in the abelian category of $\cD$-modules (index $q$ refers to the grading of the $\cA$-module $P_p\0_\cA\mathfrak{Q}\in\tt DG\cD M$). The horizontal differential $$d^h_{pq}:(P_{p}\0_\cA\mathfrak{Q})_q\to (P_{p-1}\0_\cA\mathfrak{Q})_q$$ is the $\cD$-linear map $(\bar\zp_p)_q$ given by the chain map $$\bar\zp_p:=\tilde\zp_{p}\0\id_{\mathfrak{Q}}\;.$$ The map $d^h$ actually squares to zero due to the exactness of \eqref{SequenceB}. The vertical differential $$d^v_{pq}:(P_{p}\0_\cA\mathfrak{Q})_{q}\to (P_{p}\0_\cA\mathfrak{Q})_{q-1}$$ is the $\cD$-linear map $(d_p^\0)_q$ given by the differential $$d_p^\0:=d_{P_p}\0\id+\id\0\, d_{\mathfrak{Q}}$$ of the $\cA$-module $P_p\0_\cA\mathfrak{Q}\in\tt DG\cD M\,.$ This way we actually get a double complex, i.e., $$d^v_{p-1,q}d^h_{pq}=d^h_{p,q-1}d^v_{pq}\;,$$ as $\tilde\zp_p$ is a chain map.\medskip

The second page of the horizontal spectral sequence is
$${}^h\!E^2_{pq}=H_p^h(H^v_q((P_\bullet\0_\cA\mathfrak{Q})_\bullet))=\faktor{H^v_q((P_p\0_\cA\mathfrak{Q})_\bullet)\,\cap\,\ker d^h_\sharp}{d^h_\sharp\,H^v_q((P_{p+1}\0_\cA\mathfrak{Q})_\bullet)}\;,$$ where $d^h_\sharp$ is the differential induced in homology by $d^h\,.$ Since $d^v$ is $d^\0$ and $P_i$ is a cofibrant $\cA$-module whose homology is a free $H(\cA)$-module, there is a $\tt\cD M$-isomorphism $$H^v_q((P_{i}\0_\cA\mathfrak{Q})_\bullet)\cong (H(P_i)\0_{H(\cA)}H(\mathfrak{Q}))_q\;,$$ in view of Proposition \ref{lemma2}. We thus get $${}^h\! E_{pq}^2=H^h_p((H(P_\bullet)\0_{H(\cA)}H(\mathfrak{Q}))_q)=\op{Tor}_{p}^{H(\cA)}(H(M),H(N))_q\;,$$ as $d^h_\sharp=\tilde\zp_\sharp\0\id$ when read through the previous isomorphism and as $H(X)\cong H(M)$ and $H(\mathfrak{Q})\cong H(N)\,.$\medskip

The second sheet of the vertical spectral sequence is
$${}^v\! E_{pq}^2= H^v_p(H^h_q((P_\bullet\0_\cA\mathfrak{Q})_\bullet))=\faktor{H^h_q((P_\bullet\0_\cA\mathfrak{Q})_p)\,\cap\,\ker d^v_\sharp}{d^v_\sharp\,H^h_q((P_\bullet\0_\cA\mathfrak{Q})_{p+1})}\;.$$ The horizontal homology $H^h_q((P_\bullet\0_\cA\mathfrak{Q})_\ell)$ is the homology of the differential $d^h=\tilde\zp\0\id\,,$ so that in its computation the differential $d^v=d_P\0\id +\id\0\, d_{\mathfrak{Q}}$ is irrelevant, only the graded $\cD$-module structure of $$P_k\0_\cA\mathfrak{Q}=P_k\0_\cA(\cA\0_\cO V)\in{\tt Mod}_{\tt DG\cD M}(\cA)\subset{\tt G\cD M}$$ matters. We can therefore ignore the lowering differential of $\mathfrak{Q}=\cA\0_\cO V$ and even replace it by the differential $d_\cA\0\id\,,$ thus viewing $V$ as the corresponding direct sum $S$ of spheres (this direct sum is a differential graded $\cD$-module with vanishing differential). By Proposition \ref{tri} we now get $$P_k\0_\cA\mathfrak{Q}\cong P_k\0_\cO S\in{\tt G\cD M}\;.$$ Since $\tilde\zp_k$ is $\cA$-linear, we get an isomorphism $$((P_\bullet\0_\cA\mathfrak{Q})_\ell,(\tilde\zp_\bullet\0\,\id_{\mathfrak{Q}})_\ell)\cong ((P_\bullet\0_\cO S)_\ell,(\tilde\zp_\bullet\0\,\id_S)_\ell)$$ of chain complexes of $\cD$-modules, so that $$H^h_q((P_\bullet\0_\cA\mathfrak{Q})_\ell)\cong H_q((P_\bullet\0_\cO S)_\ell)=\bigoplus_{r+s=\ell}H_q(P_{\bullet,r}\0_\cO S(s))$$ in $\tt \cD M\,.$ The complex $$P_{\bullet,r}\0_\cO S(s)=\bigoplus_kP_{k,r}\0_\cO S(s)$$ with differential $(\tilde\zp_\bullet)_r\0\id_{S(s)}$ is the tensor product of the complex $P_{\bullet,r}$ with differential $(\tilde\zp_\bullet)_r$ and the complex $S(s)$ concentrated in degree zero with zero differential. It suffices to compute its homology as homology of a tensor product of chain complexes of $\cO$-modules. We already mentioned earlier that $S(s)\,,$ $d(S(s)) = 0\,,$ $H_0(S(s))=S(s)$ and $H_{\zs>0}(S(s)) = 0$ are all $\cO$-flat. From Künneth's theorem it therefore follows that
$$H^h_q((P_\bullet\0_\cA\mathfrak{Q})_\bullet)\cong\bigoplus_\ell\bigoplus_{r+s=\ell}\bigoplus_{\zr+\zs=q}H_\zr(P_{\bullet,r})\0_\cO H_\zs(S(s))=H_q(P_\bullet)\0_\cO S\cong$$ $$H_q(P_\bullet)\0_\cA(\cA\0_\cO S)\cong H_q(P_\bullet)\0_\cA(\cA\0_\cO V)=H_q(P_\bullet)\0_\cA\mathfrak{Q}\;,$$ where the isomorphisms are isomorphisms of graded $\cD$-modules. However, the differential on the {\small LHS} is that induced by $d^v=d_P\0\id+\id\0\, d_{\mathfrak{Q}}$ and the same is true for the differential on the {\small RHS}. In other words, the third isomorphism resets the correct vertical differential, so that the {\small LHS} and the {\small RHS} are isomorphic differential graded $\cD$-modules and have therefore isomorphic vertical homologies. Since $P_\bullet$ is a resolution of $X$ in the abelian category $\tt Mod(\cA)\,,$ we have $H_{q>0}(P_\bullet)=0$ and $H_0(P_\bullet)\cong X=QM$ in ${\tt Mod}(\cA)$. Hence $${}^v\!E^2_{pq}=\{0\}\quad(q>0)$$ and $${}^v\!E^2_{p0}\cong H_p(QM\0_\cA\mathfrak{Q})\cong H_p(M\0_\cA^{\mathbb{L}}N)\;,$$ due to Equation \eqref{DerF6} of Theorem \ref{DerTenCRS}. The vertical spectral sequence thus collapses at its second page and $$H_n(\op{Tot}(C))\cong {}^v\!E^2_{n0}\cong H_n(M\0_\cA^{\mathbb{L}}N)\;$$ (see text below Equation \eqref{convergence}), so that $${}^h\!E^2_{pq}=\op{Tor}_p^{H(\cA)}(H(M),H(N))_q\Rightarrow H_{p+q}(M\0_\cA^{\mathbb{L}}N)\;,$$ as announced.

\begin{rem} \emph{Edge homomorphisms of the first quadrant double complexes' spectral sequences 
$$ G_p^v H_p(\op{Tot}C_{\bullet\bullet})\to H^v_pH^h_0(C_{\bullet\bullet})\hspace{5pt}\text{ and }\hspace{5pt}
H^h_0H^v_p(C_{\bullet\bullet})\to G_0^hH_p\op{Tot}(C_{\bullet\bullet})$$
in any Abelian category $\tt A$ are natural. More precisely, they are the $C_{\bullet\bullet}$-components of natural transformations $$G_p^v H_p\circ\op{Tot}\Rightarrow H^v_pH^h_0\quad\text{and}\quad H^h_0H^v_p\Rightarrow G_0^hH_p\circ \op{Tot}$$ between functors from the category of first quadrant double complexes in $\tt A$ to the category $\tt A$. From here, it can be shown that the edge homomorphisms of the Tor spectral sequence 
$$(H M\otimes_{H\cA}H N)_q={}^h\!E^2_{0q}\to H_q(M\otimes^\mathbb{L}_\cA N)$$
are the $N$-components of a natural transformation $$(H M\otimes_{H\cA}H -)_q\Rightarrow H_q(M\otimes^\mathbb{L}_\cA -)\;.$$}\end{rem}

\section{Exactness of the derived tensor product}\label{fm1}

Given $M\in\tt Mod(\cA)$, limit-preserving property of the tensor product $-\otimes_A M$  is not invariant under weak equivalences, and should be replaced by a homotopy invariant notion. In \cite{TV08}, the property that the derived tensor product $-\otimes_A ^\mathbb{L}M$ preserves homotopy limits is recognized as the correct notion of flatness. The rest of the paper is dedicated to the implementation of this property in our homotopical $\cD$-algebraic framework. It is worth noting that if $M$ is flat in the later derived sense, Tor spectral sequence collapses on the second page, yielding a natural isomorphism 
$$(H M\otimes_{H\cA}H N)_q\cong H_q(M\otimes^\mathbb{L}_\cA N).$$
This is the content of the Corollary \ref{last_cor}.
\subsection{Homotopy pullbacks and flat modules}

Denote by $S=\{\bullet \rightarrow \bullet \leftarrow \bullet\}$ the small category with three objects $\{c, d, b\}$, and two non-identity morphisms  $b\to d$ and $c\to d$. Given any model category $\tt M$, the functor category ${\tt Fun}(S,M)$ -- the category of cospan diagrams in $M$ -- has three model structures with objectwise week equivalences, in which fibrant objects are objectwise fibrant diagrams for which respectively $b\to d$,  $c\to d$, and  both $b\to d$,  $c\to d$ are fibrations \cite{Models}. The limit (pullback) functor is right Quillen for all three model structures, and the corresponding right derived functors, i.e. homotopy pullbacks, agree. Moreover, given a commutative square
\begin{center}
\begin{tikzcd}
A\arrow[r]\arrow[d]&B\arrow[d]\\
C\arrow[r]&D
\end{tikzcd}
\end{center}
in $\tt M$, weather the universal map $A\to \op{Lim}(R(B\to D\leftarrow C))$ from $A$ to the homotopy pullback of the cospan $B\to D\leftarrow C$ is a weak equivalence is independent of the chosen model structure or the choosen fibrant replacement $R$ in ${\tt Fun}(S,M)$. Commutative squares for which those maps are indeed weak equivalences are referred to as model squares, and denoted by $ABCD$. The vertex $A$ is called a generalzed representative of the homotopy pullback $B\times_D^h C$. In the literature, model squares are also called homotopy fiber squares \cite{Hir}, and homotopy pullback squares \cite{JL}. 

Generally, we say that a functor $F:\tt M\to N$ preserves model squares if the $F$-image of any model square in $\tt M$
is a model square in $\tt N$.
 If there is a natural weak equivalence $\eta:F\Rightarrow F'$ between two functors from $\tt M$ to $\tt N$, then by \cite[Proposition 2]{HomotopyfiberSeq} $F$ preserves model squares if and only if $F'$ does.

We now turn our attention from general functors to the derived tensor product.
To begin with, there is a certain ambiguity in its definition. Denoting by $Q$ the cofibrant F-replacement functor, and by $C$ a local cofibrant replacement, the derived tensor product
$$-\otimes^\mathbb{L}_{\cA} M:\op{Ho}({\tt Mod}({\cA}))\to\op{Ho}({\tt Mod}({\cA}))$$
can be equivalently defined as the derived functor of any of the three functors
\begin{equation}\label{dertp} Q-\otimes_{\cA} CM:N\mapsto QN\otimes_{\cA} CM,\hspace{5pt}
Q-\otimes_{\cA} M:N\mapsto QN\otimes_{\cA} M,\hspace{5pt}
-\otimes_{\cA} CM:N\mapsto N\otimes_{\cA} CM,\end{equation}
which all respect weak equivalences \ref{TPPWEQC}.

Any one of them preserves model squares if any other does, due to the natural weak equivalences
\begin{equation*}
\begin{split}
&  q\otimes_{\cA}CM:Q-\otimes_{\cA} CM\Rightarrow -\otimes_{\cA} CM,\\
&Q-\otimes c_M:Q-\otimes_{\cA} CM\Rightarrow Q-\otimes_{\cA} M.\\
\end{split}
\end{equation*}

Similarly, whether the functors in question preserve model squares is independent of the choosen cofibrant replacements. Indeed, if $Q_1, Q_2$ are two different cofibrant $F$-replacement functors, and $C$, $C'$ two different local cofibrant replacements, one gets following chains of natural weak equivalences:
\begin{equation*}
\begin{split}
&-\otimes_{\cA} CM \Leftarrow Q_1-\otimes_{\cA} CM \Rightarrow Q_1-\otimes_{\cA} M \Leftarrow Q_1-\otimes_{\cA} C'M\Rightarrow -\otimes_{\cA} C'M,\\
&-\otimes_{\cA} CM\Leftarrow Q_2-\otimes_{\cA} CM\Rightarrow Q_2-\otimes_{\cA} M\Leftarrow Q_2-\otimes_{\cA} C'M.\\
\end{split}
\end{equation*}

Finally, we are prepared to define flat modules.
\begin{defi}\label{FlatDefi}
The derived tensor product $-\otimes_{\cA}^\mathbb{L} M$ is said to \emph{preserve model squares} if any, and hence all the functors in \ref{dertp} preserve model squares. In this situation, we say that $M$ is a \emph{flat} $\cA$-module.
\end{defi}
To verify that tensoring with an $\cA$-module preserves all model squares in not a simple task. Luckily, for the derived tensor product to preserve model squares it suffices that it preserves homotopy fiber sequences. For this, we first recall the later concept, as defined in \cite{HomotopyfiberSeq}, together with the related notions which are relevant for this paper: mapping fiber, loop object (shift), and the connecting homomorphism. Definitions slightly deviate from the standard references \cite{Quill}, \cite{Ho99}.

\subsection{Homotopy fiber sequences and flat modules}

For us, a homotopy fiber sequence in a pointed model category $\tt M$, denoted by $A\to B\to D$, is a model square
\begin{equation}\label{hofibseq}
\begin{tikzcd}
A\arrow[r]\arrow[d]&B\arrow[d]\\
C\arrow[r]&D,
\end{tikzcd}
\end{equation}
in which the unique map from $C$ to the terminal/initial object is a weak equivalence. Canonical examples are fibrations together with their kernels, and the pullbacks of arbitrary maps $f:B\to D$ along the based path object (dual cone) of $D$ (\cite[Sections 4, 6]{HomotopyfiberSeq}). 
%
%
%

\medskip
To give an elegant description of the based path object it is useful to consider unbounded chain complexes as well.
A chain complex $X$ of $\cD$-modules which is not necessarily concentrated in non-negative degrees  will be called an unbounded $\cA$--module if it is equipped with a morphism $\cA\otimes X\to X$ in the category of unbounded dg $\cD$-modules, such that the usual associativity and unitality diagrams commute. The category of unbounded $\cA$--modules will be denoted by $\tt Mod (\cA)^{ub}$.  $\tt Mod (\cA)$ is its full subcategory. Given $X\in\tt Mod (\cA)^{ub}$, its good truncation, denoted by  $X^{\geq 0}$, is its subcomplex
$$\ldots\xrightarrow{d_3}D_2\xrightarrow{d_2}D_1  \xrightarrow{d_1}\op{Ker}(d_0) .$$
Being closed under the $\cA$-action, $X^{\geq 0}$ is an object in $\tt Mod(\cA)$.
It is straightforwardly verified that the obtained functor $(-)^{\geq 0}:\tt Mod (\cA)^{ub}\to Mod (\cA)$ is right adjoint to the inclusion
$\imath:\tt Mod (\cA)\to Mod (\cA)^{ub}$.

Notice that the truncation respects quasi-isomorphisms, and sends degreewise surjective maps to fibrations in $\tt Mod (\cA)$ -- the maps surjective in strictly positive degrees.

Given $D\in \tt Mod (\cA)^{ub}$,  its unbounded based path object $\op{Path}_0^{ub} D$ is its $(-1)$ -- shifted cone, explicitly,
\begin{equation*}
\begin{split} &(\op{Path}_0^{ub} D)_n= D_n\oplus D_{n+1},\\
&d=\begin{pmatrix}
d_D & 0 \\
-\op{id}_D & -d_D 
\end{pmatrix}.
\end{split}
\end{equation*}
For the explicit description of the $\cA$-action see the Remark \ref{devet}.

Denote by $\pi_D^{ub}$ the degreewise projection 
$$\op{Path}_0^{ub} D\to D,\hspace{10pt}D_n\oplus D_{n+1}\to D_n$$ 
Given a morphism $f:B\to D$ in $\tt Mod (\cA)^{ub}$, the pullback
\begin{equation}\label{puppe}
\begin{tikzcd}
{K_f^{ub}}\arrow[r,"\pi_f^{ub}"]\arrow[d,"p_f^{ub}"]&B\arrow[d,"f"]\\
\op{Path}_0^{ub} D\arrow[r,"\pi_D^{ub}"]&D,
\end{tikzcd}
\end{equation}
will be called the unbounded fiber of $f$. It equals the $(-1)$ shifted mapping cone.
 The $(-1)$-shift of $D$, 
$$D[-1]_n=D_{n+1},\hspace{10pt} d_n^{D[-1]}=-d^D_{n+1}$$
can be seen as the kernel of $\pi_D^{ub}:\op{Path}_0^{ub} D\to D$, and the connecting homomorphism $\delta^{ub}_f:D[-1]\to {K_f}^{ub}$ is the pullback's universal map
\begin{center}
\begin{equation}\label{connecting}
\begin{tikzcd}
D[-1]
\arrow[bend left]{drr}{0}\arrow[bend right,swap]{ddr}{\op{Ker}}\arrow[dashed]{dr}{\delta^{ub}_f}
& & \\
& K_f^{ub} \arrow[r,"\pi^{ub}_f"] \arrow[d,"p_f^{ub}"] & B \arrow[d,"f"] \\
& \op{Path}_0^{ub} D \arrow[r,"\pi^{ub}_D"] & D.
\end{tikzcd}
\end{equation}
\end{center}

Assume now $f:B\to D$ is a morphism in $\tt Mod(\cA).$ By a slight abuse of notation, we denote the map $\imath f$ in $\tt Mod(\cA)^{ub}$  also by $f$. Let 
\begin{center}
\begin{equation}\label{connecting_tr}
\begin{tikzcd}
\Omega D
\arrow[bend left]{drr}{0}\arrow[bend right,swap]{ddr}{\op{Ker}}\arrow[dashed]{dr}{\delta_f}
& & \\
& K_f \arrow[r,"\pi_f"] \arrow[d,"p_f"] & B \arrow[d,"f"] \\
& \op{Path}_0 D \arrow[r,"\pi_D"] & D.
\end{tikzcd}
\end{equation}
\end{center}
be the image of Diagram (\ref{connecting}) under truncation.

As truncation is right adjoint, the commutative square $K_f, B, \op{Path_0} D, D$ is a pullback. As $\op{Path_0} D$ is acyclic, and $\pi_D$ is a fibration, it is even a homotopy fiber sequence. Finally, homotopy fiber sequences induce long exact sequences in homology, in this example
$$\cdots\xrightarrow{H_1 f} H_1(D)\xrightarrow{H_0\delta_f}H_0({K_f})\xrightarrow{H_0 p_f} H_0(B)\xrightarrow{H_0 f} H_0(D).$$

Looping functor $\Omega$, defined as the composition of $(-1)$--shift and truncation has the left adjoint suspension functor $\Sigma$, which equals the $1$--shift.
Explicitly, it is defined on objects as
$$\Sigma D_n= D_{n-1},\hspace{5pt}d^{\Sigma D}_n=-d^D_{n-1},$$
and on morphisms as
$$\Sigma f_n= f_{n-1}.$$
Since $\Omega$ presereves fibrations and weak equivalences, $\Sigma\dashv \Omega$ is even a Quillen adjunction. Unlike with unbounded chain complexes, where $1$--shift is the two sided inverse of $(-1)$--shift, $\Sigma$ is a only a right inverse to $\Omega$.\bigskip

\begin{rem}\label{devet} \em{ There is a subtlety involved when determining the correct $\cA$-action on $\op{Path}_0^{ub} D$ (which induces the $\cA$-action on both $\op{Path}_0D$ and $\Omega D$).  $\cA$ being graded commutative, the notions of left and right module structures on $D$ are equivalent, and related by
$m\triangleright a=(-1)^{am}a\triangleleft m$
(Remark \ref{sest}). A reasonable way to define the $\cA$-action on $\op{Path}_0^{ub} D$ is from the $\cA$-action on $D$, in the decomposition $\op{Path}_0^{ub} D=D\oplus D[-1]$. A priori, there is a choice: either one takes the left induced action (coming from the term-wise left $\cA$-action), or the right induced action. The two will agree on the first summand, but differ by a sign on the second. However, only the right induced action is compatible with the differential. Concretely, for $m\in (\op{Path}_0 D)_n$, and $a\in \cA$, denoting by $\cdot$ the $\cA$-action on $D$, $\cA$-action on  $(\op{Path}_0 D)_n$, denoted by $\triangleleft$, is given by
$$ a\triangleleft m=\begin{cases} a\cdot m, & m\in D_n\\  (-1)^a a\cdot m, &m\in D_{n+1}.  \end{cases}$$

Consequently, $\cA$-action on the $(-1)$--shift $D[-1]$, and its truncation $\Omega D$ is given by $a\triangleleft m=(-1)^aa\cdot m$. To assure that unit and counit of $\Sigma\dashv\Omega$ adjunction respect the $\cA$-action, the same formula defines the $\cA$-action on the $(1)$--shift  $\Sigma D$.}\end{rem}

With all the relevant notions illuminated, we move on to showing that  $-\otimes^\mathbb{L}_A M$ preserves all model squares if it preserves all homotopy fiber sequences. Again, there is a certain ambiguity in what it means for the derived tensor product to preserve homotopy fiber sequences. 

In general, we say that a functor $F:\tt M\to N$ preserves homotopy fiber sequences if the $F$-image of any homotopy fiber sequence in $\tt M$ (diagram \ref{hofibseq})
is a homotopy fiber sequence in $\tt N$.  If there is a natural weak equivalence $\eta:F\Rightarrow F'$ between two functors from $\tt M$ to $\tt N$, then by \cite[Proposition 2]{HomotopyfiberSeq} $F$ preserves homotopy fiber sequences if and only if $F'$ does.  Thus, the following definition makes sense:
\begin{defi}\label{SeqPul}
The derived tensor product $-\otimes_{\cA}^\mathbb{L} M$ is said to \emph{preserve homotopy fiber sequences} if any, and hence all the functors in \ref{dertp} preserve homotopy fiber sequences.
\end{defi}
Similarly, as detailed in the previous section, any of the functors \ref{dertp} representing $-\otimes^\mathbb{L}_A M$ preserves model squares whenever any other does. Consequently, to show that $-\otimes^\mathbb{L}_A M$ preserves all model squares if it preserves all homotopy fiber sequence, it suffices to prove the statement for the functor $-\otimes_{\cA} CM$, or equivalently, that for any cofibrant object $M\in {\tt Mod}({\cA})$, the functor $-\otimes_{\cA} M$ preserves all model squares if it preserves all homotopy fiber sequences.
\begin{prop}
A functor $F:{\tt Mod}({\cA})\to {\tt Mod}({\cA})$ which preserves weak equivalences  preserves all model squares if it preserves all homotopy fiber sequences.
\end{prop}
\begin{proof}
Let
\begin{center}
\begin{equation}\label{pbl}
\begin{tikzcd}
A\arrow[r]\arrow[d]&B\arrow[d,"f"]\\
C\arrow[r,"g"]& D
\end{tikzcd}
\end{equation}
\end{center}
be a model square in ${\tt Mod}({\cA})$.
Given any $\tt DG\cD M$--mmorphism $\phi:X\to Y$, the map $\Omega \phi$ fits in a model square
\begin{center}
\begin{equation}\label{pbr}
\begin{tikzcd}
\Omega X\arrow[r]\arrow[d,"\Omega \phi"]&\op{Path}_0 X \arrow[d]\\
\Omega Y\arrow[r]&K_\phi
\end{tikzcd}
\end{equation}
\end{center}
(see the diagram 49 in \cite{HomotopyfiberSeq}). $\Sigma$ being the right inverse of $\Omega$, the model squares \ref{pbl} and \ref{pbr} (for $\phi=\Sigma f$) can be pasted into
\begin{center}
\begin{tikzcd}
A\arrow[r]\arrow[d]&B\arrow[d,"f"]\arrow[r]&\op{Path}_0 \Sigma B \arrow[d]\\
C\arrow[r,"g"]&D\arrow{r}&K_{\Sigma f}.
\end{tikzcd}
\end{center}
By the pasting law \cite[Proposition 8]{Models}, the total square is a model square as well. It is in fact a homotopy fiber sequence, as $\op{Path}_0 \Sigma B$ is acyclic.
Applying the functor $F$ which by the assumption preserves homotopy fiber sequences, we get a diagram
\begin{center}
\begin{tikzcd}
FA\arrow[r]\arrow[d]&FB\arrow[d,"Ff"]\arrow[r]&F\op{Path}_0 \Sigma B \arrow[d]\\
FC\arrow[r,"Fg"]&FD\arrow{r}&FK_{\Sigma f}
\end{tikzcd}
\end{center}
whose right-hand-side square and total square are homotopy fiber sequences, especially model squares. By the pasting law we conclude that the left-hand-side square is a model square as well.
\end{proof}

\begin{cor}\label{FlatSeq}
$M\in \tt Mod(\cA)$ is flat if and only if $-\otimes^\mathbb{L}_{\cA} M$ preserves homotopy fiber sequences.\end{cor}
\section{A characterization of flat modules}\label{fm2}
We begin with the definition of strongly flat modules.
\begin{defi}\label{Strong}\begin{enumerate}
              \item Let $\cA\in\tt DG\cD A$. A \emph{module} $M\in\tt Mod(\cA)$ is {\bf strong}, if the natural $\,{\tt Mod_{G\cD M}}(H_0(\cA))$-morphism \be\label{NatMorphTens}\phi_{\bullet,\,\cA, M}:H_\bullet(\cA)\0_{H_0(\cA)}H_0(M)\to H_\bullet(M), \hspace{5pt}[a]\otimes [m_0]\mapsto [a\cdot m_0]\ee is an isomorphism, i.e., if all $\,{\tt Mod_{\cD M}}(H_0(\cA))$-morphisms $\phi_{k,\,\cA,M}$ are isomorphisms.

              \item $M\in\tt Mod(\cA)$ is {\bf strongly flat} if it is strong, and if $H_0(M)$ is a flat $H_0(A)$--module in the classical sense.
            \end{enumerate}
\end{defi}
Before moving forward with the proof that strongly flat modules are equivalently flat, we prove that the isomorphism (\ref{NatMorphTens}) is well defined.

 Let $\cA\in \cD\tt A$. In order to define a ${\tt Mod_{\cD M}}(\cA)$-morphism (an $\cA$- and $\cD$-linear map) $\phi:M'\0_\cA M''\to M$, one starts defining an $\cO$-bilinear map $\varphi:M'\times M''\to M$, hence, an $\cO$-linear map $\varphi:M'\0_\cO M''\to M$. One then checks that $\varphi$ is a $\tt \cD M$-morphism, i.e., is $\cD$-linear, or, equivalently, is linear for the action $\nabla_\theta$ of vector fields $\theta$. Recall that, by definition,
$$
\nabla_\theta(m'\0 m'')=(\nabla_\zy m')\0 m'' + m'\0(\nabla_\zy m'')\;.
$$

Denote by $\zk:M'\0_\cO M''\to M'\0_\cA M''$ coequalizer's universal map.  If, for $a\in\cA$,
$$
\varphi((a\cdot m')\0 m'')=\varphi(m'\0(a\cdot m''))\;,
$$
it follows from the universal property of a coequalizer, that there is a unique $\cD$-linear map $\phi:M'\0_\cA M''\to M$, such that $\phi\circ \zk=\varphi$. Finally, if $\varphi$ is $\cA$-bilinear on $M'\times M''$, then $\phi$ is $\cA$-linear on $M'\0_\cA M''$. Indeed,
$$
\varphi(m',m'')=\varphi(m'\0 m'')=\phi(\zk(m'\0 m''))=\phi(m'\0_\cA m'')\;.
$$

We now come back to Definition \ref{Strong}. The tensor product in \eqref{NatMorphTens} makes sense since $H_0(\cA)\in\cD A$ and $H_k(\cA), H_0(M)\in{\tt Mod_{\cD M}}(H_0(\cA))$. In order to define the morphism $\phi_{k,\,\cA,M}$ (we will write $\phi$), we apply the just detailed method to the preceding $\cD$-algebra and modules in $\tt\cD M$ over that algebra. Let $\varphi$ be the map $$\varphi:H_k(\cA)\times H_0(M)\ni ([a_k],[m_0])\mapsto [\zn(a_k \0 m_0)]\in H_k(M)\;,$$ where the $\tt DG\cD M$-morphism $\zn:\cA\0_\cO M\to M$ is the action of $\cA$ on $M$. It is easy to check that $\varphi$ is well-defined on homology classes. In view of the $\tt\cD A$-morphism $$\cO\ni f\mapsto f\cdot[1_\cA]=[f\cdot 1_\cA]\in H_0(\cA)\;,$$ if $\varphi$ is $H_0(\cA)$-bilinear, it is in particular $\cO$-bilinear. As for $H_0(\cA)$-bilinearity, taking into account that the action $\ast$ of $H_0(\cA)$ on $H_k(\cA)$ is induced by the multiplication $\ast$ in $\cA$, and that the action $\zn$ of $H_0(\cA)$ on $H_k(M)$ is induced by the action $\zn$ of $\cA$ on $M$, we get
$$
\varphi([a_0]\ast[a_k],[m_0])=\varphi([a_0\ast a_k],[m_0])=[\zn((a_0\ast a_k)\0 m_0)]\;,
$$
$$
\zn([a_0]\0\varphi([a_k],[m_0]))=\zn([a_0]\0[\zn(a_k\0 m_0)])=[\zn(a_0\0\zn(a_k\0 m_0))]=[\zn((a_0\ast a_k)\0 m_0)]\;,
$$
and
$$
\varphi([a_k],\zn([a_0]\0[m_0]))=\varphi([a_k],[\zn(a_0\0 m_0)])=[\zn(a_k\0\zn(a_0\0 m_0))]=[\zn((a_k\ast a_0)\0 m_0)]\;,
$$
so that $\varphi$ is actually $H_0(\cA)$-bilinear and thus $\cO$-bilinear. We now check linearity of $$\varphi:H_k(\cA)\0_\cO H_0(M)\to H_k(M)$$ with respect to the action $\nabla_\zy$ by vector fields $\zy$. As the $\cD$-action on homology is induced by the $\cD$-action on the underlying complex, we have
$$
\zvf(\n_\zy([a_k]\0 [m_0]))=\zvf((\n_\zy[a_k])\0 [m_0])+[a_k]\0(\n_\zy[m_0]))=
$$
$$
\zvf([\n_\zy a_k]\0 [m_0])+\zvf([a_k]\0[\n_\zy m_0])=[\zn((\n_\zy a_k)\0 m_0)]+[\zn(a_k\0(\n_\zy m_0))]=
$$
$$
[\zn(\n_\zy(a_k\0 m_0))]=\n_\zy[\zn(a_k\0 m_0)]=\n_\zy\,\zvf([a_k]\0[m_0])\;.
$$
In view of what has been said above, it follows that $\zvf$ induces a unique ${\tt Mod_{\cD M}}(H_0(\cA))$-morphism $\zf:H_k(\cA)\0_{H_0(\cA)}H_0(M)\to H_k(M)$.

\subsection{Stronly flat modules are flat}\label{flat1}
Notice that both flatness and strong flatness are invariant under the isomorphisms in the homotopy category. This allows us to assume without the loss of generality that $M$ is cofibrant, so that
 $-\otimes_{\cA} M$ represents the derived tensor product $-\otimes_{\cA}^\mathbb{L} M$.


Assume that $M$ is a strongly flat $\cA$-module. To prove that $M$ is flat, it suffices to show that the functor $-\otimes_{\cA} M$ preserves homotopy fiber sequences. By  \cite[Corollary 7]{HomotopyfiberSeq}, we further reduce to showing that for any map $f:B\to D$ of $\cA$-modules, the image of the homotopy fiber sequence 
\begin{equation}
\begin{tikzcd}
{K_f}\arrow[r,"\pi_f"]\arrow[d,"p_f"]&B\arrow[d,"f"]\\
\op{Path}_0D\arrow[r,"\pi_D"]&D,
\end{tikzcd}
\end{equation}
under the functor $-\otimes_{\cA} M$ is still a homotopy fiber sequence.

Let $f:B\to D$ be a morphism in $\tt Mod(\cA)$. Unbounded $\cA$-modules $D[\pm1]\otimes M$ and $(D\otimes M)[\pm1]$ are equal, since the equality of graded modules
$$(D[\pm1]\otimes M)_n=\bigoplus_{k+l=n}(D_{k\pm1}\otimes M_l)=\bigoplus_{k+l=n\pm1}(D_{k}\otimes M_l)=(D\otimes M)[\pm1]_n$$
respects both the differential and the $\cA$-action. As right $\cA$--actions on $D[\pm1]$ and $D$ coincide (Remark \ref{devet}), 
\begin{equation}\label{shifts}D[\pm1]\otimes_\cA M=(D\otimes_\cA M)[\pm1].\end{equation}
On the other side, the tensor products $M\otimes_\cA D[\pm1]$ and $(M\otimes_\cA D)[\pm1]$ not equal, but can still be identified via the isomorphism $m\otimes d\mapsto (-1)^{m}m\otimes d$.

The isomorphism of graded $\cA$-modules
\begin{equation*}
(\op{Path}_0^{ub} D)\otimes_{\cA} M=(D\oplus D[-1])\otimes_\cA M\to D\otimes_\cA M \oplus D[-1]\otimes_\cA M=\op{Path}_0^{ub}(D\otimes_{\cA} M)
\end{equation*}
also respects the differential. The obtained isomorphism in $\op{Mod}^{ub}(\cA)$ is $$\phi^{ub}:(\op{Path}_0^{ub} D)\otimes_{\cA} M\xrightarrow{\cong} \op{Path}_0^{ub}(D\otimes_{\cA} M), \hspace{5pt}(d+d^{-1})\otimes m\mapsto d\otimes m+d^{-1}\otimes m,$$
for $d\in D$, $d^{-1}\in D[-1]$, $m\in M.$

Denote by $K^{\otimes,ub}$ the pullback
$$(B\otimes_{\cA} M)\prod_{D\otimes_{\cA} M} \op{Path}^{ub}_0(D\otimes_{\cA} M),$$ and by $K^\otimes$ its good truncation. Denote by $\epsilon$ the counit of the $\imath\vdash (-)^{\geq 0}$ adjunction, explicitly the natural map
\begin{center}
\begin{equation}\label{big}
\begin{tikzcd}
\ldots\arrow[r,"d_3"]& D_2\arrow[r,"d_2"]\arrow[d,equal]&D_1\arrow[d,equal]\arrow[r,"d_1"]&\op{Ker}(d_0)\arrow[r]\arrow[d,hook]&0\arrow[r]\arrow[d,]&\ldots \\
\ldots\arrow[r,"d_3"]& D_{2}\arrow[r,"d_{2}"]&D_1\arrow[r,"d_1"] &D_0\arrow[r,"d_0"]&D_{-1}\arrow[r,"d_{-1}"]&\ldots
\end{tikzcd}
\end{equation}
\end{center}
As $(\op{Path}_0 D)\otimes_{\cA} M$ is non-negatively graded, the composition
$$(\op{Path}_0 D)\otimes_{\cA} M\xrightarrow{\epsilon\otimes M}(\op{Path}_0^{ub} D)\otimes_{\cA} M\xrightarrow{\phi^{ub}} \op{Path}_0^{ub} (D\otimes_{\cA} M)$$
factors through $\op{Path}_0 (D\otimes_{\cA} M).$ Let
$$\phi:(\op{Path}_0 D)\otimes_{\cA} M\rightarrow \op{Path}_0(D\otimes_{\cA} M)$$
be the resulting map.

The following commutative diagram, in which the dashed arrows are universal maps of the pullback, is central. For notation, see Diagrams (\ref{connecting}) and (\ref{connecting_tr}).
\begin{center}
\begin{equation}\label{big}
\begin{tikzpicture}[baseline= (a).base]
\node[scale=.9] (a) at (0,0){
\begin{tikzcd}
\Omega D\otimes_{\cA} M \arrow[rrr]\arrow[dr,dashed,"\psi_\Omega"]\arrow[dd,"\delta_f\otimes M"]&&&0\arrow[dd]\\
&\Omega(D\otimes_{\cA} M)\arrow[rru]\arrow[r,"\epsilon"]\arrow[dd]\arrow[dd,"\delta_{f\otimes M}"{yshift=6.5pt}]&(D\otimes_{\cA} M)[-1]\arrow[ru]\arrow[dd,"\delta^{ub}_{f\otimes M}"{yshift=6.5pt}]&\\
K_f\otimes_{\cA} M \arrow[dd,"p_f\otimes M"]\arrow[rrr,"\pi_f\otimes M"]\arrow[dr,dashed,"\psi_K"]&&&B\otimes_{\cA} M\arrow[dd,"f\otimes M"]\\
& K^\otimes\arrow[dd]\arrow[rru,"\pi_{f\otimes_\cA M}"]\arrow[r,"\epsilon"]\arrow[dd,"p_{f\otimes M}"{yshift=6.5pt}]& K^{\otimes,ub}\arrow[dd,"p^{ub}_{f\otimes M}"{yshift=6.5pt}]\arrow[ru,swap,"\pi_{f\otimes_\cA M}^{ub}"]\arrow[dd]&\\
(\op{Path}_0 D)\otimes_{\cA} M\arrow[dr,"\phi"] \arrow[dr]\arrow[rrr,"\pi_D\otimes M"]&&& D\otimes_{\cA} M\\
& \op{Path}_0(D\otimes_{\cA} M)\arrow[rru,"\pi_{D\otimes_\cA M}"]\arrow[r,"\epsilon"]&\op{Path}^{ub}_0(D\otimes_{\cA} M)\arrow[ru,swap,"\pi_{D\otimes_\cA M}^{ub}"]&
\end{tikzcd}};
\end{tikzpicture}
\end{equation}
\end{center}

To prove that $-\otimes M$ preserves homotopy fiber sequences, it suffices to show that $\psi_K$ is a quasi-isomorphism. The composition $p_f\circ \delta_f$ is the inclusion of the kernel of $\pi_D$ (Diagram \ref{connecting_tr}). Likewise, $p^{ub}_{f\otimes M} \circ \delta^{ub}_{f\otimes M} $ is inclusion of the kernel of $\pi^{ub}_{D\otimes_\cA M}$. Finally, the composition $\epsilon\circ \psi_\Omega$ equals the natural map $$\epsilon\otimes M: \Omega D\otimes_\cA M\to D[-1]\otimes_\cA M=(D\otimes_\cA M)[-1],$$ as, in view of kernel's universal property, it is the unique map for which the square 
\begin{center}
\begin{tikzcd}
[column sep={7em}, row sep={3em}]
\Omega D\otimes_\cA M \arrow[r,"\epsilon \circ \psi_{\Omega}"]\arrow[r,swap,"=\epsilon\otimes M"]\arrow[d,"\substack{(p_f\otimes M) \circ (\delta_f\otimes M)\\=\op{ker}(\pi_D)\otimes M}"]&D[-1] \otimes_\cA M\arrow[d,"\substack{p^{ub}_{f\otimes M} \circ \delta^{ub}_{f\otimes M}\\=\op{ker}(\pi_D^{ub})\otimes M} "]\arrow[r,equal]&(D\otimes_\cA M)[-1]\arrow[d,"\op{ker}(\pi^{ub}_{D\otimes M})"]\\
\op{Path}_0 D\otimes_\cA M \arrow[rr, bend right=7,swap, "\epsilon\circ \phi"]\arrow[r,"\epsilon\otimes M"]&\op{Path}_0^{ub}D\otimes_\cA M\arrow[r,"\phi^{ub}"]&\op{Path}_0^{ub}(D\otimes_\cA M)
\end{tikzcd}
\end{center}
commutes. The left-hand side commutes due to the naturality of $\epsilon$. The commutativity of the right-hand side follows from the explicit description of involved maps. 

Strongness of $M$, together with the Proposition \ref{tri}, gives the isomorphisms
\begin{equation}\label{first} H_\bullet D \otimes_{H_0 \cA} H_0M\cong H_\bullet D\otimes_{H_\bullet \cA}(H_\bullet \cA\otimes_{H_0 \cA}H_0 M)\cong H_\bullet D\otimes_{H_\bullet \cA} H_\bullet M,\end{equation}
natural in $D$.
As $H_0M$ is a flat $H_0\cA$-module, $H_\bullet M$ is a flat graded $H_\bullet \cA$-module, and the Tor spectral sequence
$$\op{Tor}_p^{H_\bullet \cA}(H_\bullet D, H_\bullet M)_q \Rightarrow H_{p+q}(D\otimes_\cA M)$$
collapses on the second page. Consequently, the natural edge homomorphism 
\begin{equation}\label{quil}
 H_\bullet D \otimes_{ H_\bullet \cA} H_\bullet M\to H_\bullet(D\otimes_{\cA} M)\end{equation} is an isomorphism. Composing with \ref{first}, we get the isomorphism 
\begin{equation} \label{second} H_\bullet D  \otimes_{ H_0 \cA} H_0 M\to H_\bullet(D\otimes_{\cA} M),\end{equation}
natural in $D$.
Applying to $\epsilon$, we conclude that 
$$H_k(\epsilon\otimes M): H_k(\Omega D\otimes_\cA M)\to H_k(D[-1]\otimes_\cA M)$$
is an isomorphism for $k\geq 0.$ A word of caution is in order: since $\epsilon:\Omega D\to D[-1]$ is not a map of complexes concentrated in non-negative degrees, to be completely rigorous, before using the natural isomorphism (\ref{second}), one should apply the $1$--shift.

Being the good truncation of $\epsilon\otimes M$, $\psi_\Omega$ is a quasi-isomorphism, natural in $D$. In the diagram
\begin{center}\begin{tikzpicture}
\node[scale=.68] (a) at (0,0){
\begin{tikzcd}[every label/.append style = {font = \small},column sep=small]
H_{\bullet+1}(B) \otimes_{H_0(A)}H_0(M)\arrow[r]\arrow[d,"\cong"]& H_{\bullet+1}(D)\otimes_{H_0(A)}
H_0(M) \arrow[r]\arrow[d,"\cong"] & H_\bullet(K_f)\otimes_{H_0(A)}
H_0(M) \arrow[r]\arrow[d,"\cong"] &H_\bullet(B)\otimes_{H_0(A)}H_0(M) \arrow[r]\arrow[d,"\cong"]& H_\bullet(D) \otimes_{H_0(A)}H_0(M)\arrow[d,"\cong"] \\
H_\bullet(\Omega B\otimes_{\cA} M) \arrow[r]\arrow[d,"H_\bullet(\psi_{\Omega})"]&
H_\bullet(\Omega D\otimes_{\cA} M) \arrow[r]\arrow[d,"H_\bullet(\psi_{\Omega})"]& H_\bullet( K_f\otimes_{\cA} M) \arrow[r]\arrow[d,"H_\bullet(\psi_K)"] &H_\bullet( B\otimes_{\cA} M) \arrow[r]\arrow[d,equal]&H_\bullet( D\otimes_{\cA} M) \arrow[d,equal] \\
H_\bullet \Omega (B\otimes_{\cA} M) \arrow[r]&
H_\bullet \Omega (D\otimes_{\cA} M) \arrow[r]& H_\bullet( K^\otimes) \arrow[r] &H_\bullet( B\otimes_{\cA} M) \arrow[r]&H_\bullet( D\otimes_{\cA} M)\\
\end{tikzcd}};\end{tikzpicture}\end{center}
the bottom line is the exact sequence in homology associated to the homotopy fiber sequence 
$$K^\otimes\to B\otimes_{\cA} M\to D\otimes_{\cA} M.$$
The top line is the long exact sequence in homology associated to the homotopy fiber sequence
$$K_f\to B\to D,$$
tensored with the  flat $H_0(A)$-module $H_0(M)$. Since $\psi_\Omega$ is a quasi-isomorphism, 5-lemma implies that $\psi_K$ is a quasi-isomorphism as well.

\subsection{Postnikov towers}
The proof that all flat modules are strognly flat uses the machinery of Postnikov towers, which is  -- for the category ${\tt Mod}({\cA})$ -- recalled in this section.

For given $M\in {\tt Mod}({\cA})$, set $M_{\leq n}$ to be the complex
\begin{equation}\label{postnikov}\ldots\rightarrow 0\rightarrow \op{Im}(d_{n+1})\hookrightarrow M_{n}\xrightarrow{d_n} M_{n-1}\xrightarrow{d_{n-1}}\ldots \xrightarrow{d_1} M_0.\end{equation}
Degreewise surjection $\phi_n:M\to M_{\leq n}$ is defined in degrees $\leq n$ as identity, and in degree $n+1$ as $d_{n+1}$. Degreewise surjections $p_n:M_{\leq n}\to M_{\leq n-1}$ are defined alike. 
$M_{\leq n}$ is an $\cA$-module, with the $\cA$-action defined as follows: 
\begin{itemize}
\item
For $a\in A$ and $m\in M_k, k\leq n$, we set $a\cdot m$ to be image of $a\cdot m\in M$ under $\phi_n$. 
\item For $dm\in \op{Im}(d_{n+1})$, $a\cdot dm$ is set to zero unless $a\in A_0$, and for $a\in A_0$ we have $a\cdot dm=d(a\cdot m)\in \op{Im}(d_{n+1})$.
\end{itemize}

Maps $p_n$ and $\phi_n$ are compatible in the sense that $ p_n\circ\phi_n=\phi_{n-1}$, and are easily checked to be morphisms of $\cA$-modules. In fact, being degreewise surjective, the maps are fibrations.
\begin{defi}
Postnikov tower of $M\in {\tt Mod}({\cA})$ is the inverse system 
$$\ldots\to M_{\leq n}\xrightarrow{p_n} M_{\leq n-1}\xrightarrow{p_{n-1}}\ldots\xrightarrow{p_2} M_{\leq 1}\xrightarrow{p_1} M_{\leq 0},$$
together with the sequence of compatible maps $\phi_n:M\to M_{\leq n}$.
\end{defi}

Indeed, so-defined Postnikov towers satisfy the defining properties of the more familiar topological concept:

\begin{prop}
\begin{enumerate}
\item $H_k(\phi_n):H_k(M)\to H_k(M_{\leq n})$ is an isomorphism for $k\leq n$;
\item $H_k(M_{\leq n})=0$ for $k>n$;
\item the homotopy fiber of the fibration $M_{\leq n}\to M_{\leq n-1}$ is $H_{n}(M)$, viewed as a complex concentrated in degree $n$.
\end{enumerate}
\end{prop}
\begin{proof}
The first two statements are obvious. For the third, it suffices to notice that the homotopy fiber in question coincides with the kernel of $p_n$ (on the nose), as the map in question is a fibration, and all objects in ${\tt Mod}({\cA})$ are fibrant.
\end{proof}

As with CW complexes, $M$ can be recovered as the homotopy limit of the inverse system $(M_{\leq n})_{n\in \mathbb{N}}$. Namely, given any model category $\tt M$, the functor category $\op{Fun}(\mathbb{N}^{\op{op}},\tt M)$ has an injective model structure whose weak equivalences and cofibrations are objectwise \cite[Theorem 5.3.1.]{Ho99}. Fibrant objects are sequences of fibrations between fibrant objects in $\tt M$. Especially, Postnikov tower $\ldots M_{\leq n}\twoheadrightarrow M_{\leq n-1}\twoheadrightarrow\ldots\twoheadrightarrow M_{\leq 0}$ is a fibrant object in $\op{Fun}(\mathbb{N}^{\op{op}},{\tt Mod}({\cA}))$, so that its homotopy limit coincides with  its projective limit on the nose, which is clearly $M$.

\subsection{Flat modules are strongly flat}\label{flat2}

 Given a morphism $f:\cA\to\cB$ in $\tt DG\cD A$, any $\cB$-module has an induced $\cA$-module
structure defined as $a\cdot m:=f(a)\cdot m$. For $\cA\in \tt DG\cD A$, both $\cA_0$, and $H_0\cA$ are commutative $\cD$-algebras, and can be viewed as differential graded $\cD$-algebras concentrated in degree zero. Due to $\tt DG\cD A$ morphisms $\cA_0\hookrightarrow\cA\twoheadrightarrow H_0\cA$, any $H_0\cA$-module is also an $\cA$--module, and any $\cA$--module is also an $\cA_0$--module. These facts will be used repeatedly throughout this subsection. For $D\in\tt Mod(\cA)$, we will denote by $\pi_0:D\to H_0(D)$ the natural map
\begin{center}
\begin{equation}\label{big}
\begin{tikzcd}
 \ldots\arrow[r,"d_3"]&D_2\arrow[d]\arrow[r,"d_2"]&D_1\arrow[d]\arrow[r,"d_1"]&D_0\arrow[d,two heads]\\
\dots\arrow[r]&0\arrow[r] &0\arrow[r]&H_0(D)
\end{tikzcd}
\end{equation}
\end{center}

Suppose that $M$ is a flat $\cA$-module. As  in  the subsection \ref{flat1}, we assume $M$ to be cofibrant. Any exact sequence 
$$0\to N\hookrightarrow  P$$
 of $H_0\cA$-modules is a homotopy fiber sequence in $\tt Mod(\cA)$, since any map between complexes concentrated in degree zero is a fibration. Since $M$ is flat,
\begin{equation}\label{fiseq}0\to N\otimes_\cA M\to P\otimes_\cA M\end{equation}
is a homotopy fiber sequence as well.
From the associated long exact sequence in homology, we conclude that
\begin{equation}\label{ses} 0\to H_0(N\otimes_{\cA} M)\rightarrow H_0(P\otimes_{\cA} M)\end{equation}
is an exact sequence of $H_0\cA$-modules.

\begin{lemma}
For $N\in{\tt Mod }(H_0\cA)$, $M\in \tt Mod(\cA)$, there exists a natural isomorphism $$\phi:H_0( N\otimes_{\cA} M)\to N\otimes_{H_0\cA} H_0M,\hspace{10pt}[n\otimes m_0]\mapsto n\otimes [m_0]$$ of $H_0\cA$--modules.
\end{lemma}
From here and (\ref{ses}), it follows that $H_0M$ is a flat $H_0 \cA$--module.
\begin{proof}
For $\cA\in \tt DG\cD A$, and $P,Q\in \tt Mod(\cA)$, the tensor product $P\otimes_\cA Q$ is the quotient of $P\otimes Q$ by the graded submodule $\cI_A$ generated by the homogeneous elements $p\cdot a\otimes q-p\otimes a\cdot q$, together with the induced differential, and the $\cA$-action induced by the left $\cA$-action on $P$. For $P,Q\in \tt Mod(\cB)$, the inclusion $\cI_A\subseteq \cI_B$ of submodules in $P\otimes Q$ induces the quotient map
$P\otimes_\cA Q\twoheadrightarrow P\otimes_\cB Q$ in $\tt Mod(\cB)$.

For $m\in M$, $a\in\cA$, and $n\in N$, $n\cdot a\otimes m-n\otimes a\cdot m$ is a homogenious element of degree zero in $M\otimes N$ if and only if both $m$ and $a$ are of degree zero. Thus, in degree zero, the projection $p:N\otimes_{\cA_0}M\twoheadrightarrow N\otimes_{\cA}M$ is isomorphism. As, $p$ is degreewise surjective, $H_0 p$ is an isomorphism of $H_0\cA$-modules. Due to right exactness of the tensor product, the universal map
\begin{equation*}
\begin{split}  H_0( N\otimes_{\cA_0}M)= &\op{coker}(N\otimes_{\cA_0}M_1 \xrightarrow{\op{id}\otimes d} N\otimes_{\cA_0}M_0)\\ 
\to & N\otimes_{\cA_0} \op{coker}(M_1\xrightarrow{d} M_0)=N\otimes_{\cA_0} H_0M,\\
&[n\otimes m_0]\mapsto n\otimes [m_0]
\end{split}
\end{equation*}
is an isomorphism.
Composing with inverse of $H_0 p$, we get another isomorphism
\begin{equation}\label{lema2}H_0(N\otimes_{\cA}M)\to N\otimes_{\cA_0} H_0M,\hspace{10pt} [n\otimes m_0]\mapsto n\otimes [m_0].\end{equation}

As the map $\cA_0\twoheadrightarrow H_0(\cA)$ is surjective, the quotient map $$N\otimes_{\cA_0} H_0M\to N\otimes_{H_0\cA} H_0M$$ is an isomorphism as well. Composing with (\ref{lema2}), we get the desired result.

\end{proof}

We now show that that $M$ is a strong $\cA$-module.
\begin{prop}
Let $M, N\in\tt Mod(\cA)$, with $M$ flat and cofibrant. There exists a natural isomorphism of graded $H_0(\cA)$-modules
\begin{equation}\label{propo}\phi_{\bullet, N, M}:H_\bullet(N)\otimes_{H_0(\cA)}H_0(M)\to H_\bullet (N\otimes_\cA M),\hspace{10pt}[n]\otimes[m_0]\mapsto[n\otimes m_0].\end{equation}
 Applying to $N=\cA$, we find that $M$ is a strong $\cA$-module.
\end{prop}
\begin{proof}
Taking $N=0$, from the long exact sequence in homology associated to the fiber sequence \ref{fiseq}, we conclude that $H_i(P\otimes_{\cA} M)=0$, for any $H_0\cA$-module $P$ and $i>0$. Thus, the map
\begin{equation}\label{base}\pi_0:P\otimes_{\cA} M\to H_0(P\otimes_{\cA} M) \end{equation}
is a weak equivalence. Since $(n)$--shifts preserve weak equivalences, for any $H_0(A)$-module $P$, and $n>0$, the map
\begin{equation} \label{shifted} P[n]\otimes_{\cA} M= (P\otimes_{\cA}M)[n]\xrightarrow{(\phi\circ \pi_0)[n]} (P\otimes_{H_0 \cA} H_0M)[n]=P[n]\otimes_{H_0 \cA} H_0M\end{equation}
which is in degree $n$ given by  $$p\otimes m_0\mapsto p\otimes [m_0]$$
is a weak equivalence as well.

Notice that for any $n\in\mathbb{N},$ and $i< n$, 
\begin{equation}\label{bezveze}H_i(\phi_n\otimes \op{id}_M ):H_i(N\otimes_{\cA} M)\to H_i(N_{\leq n}\otimes_{\cA} M)\end{equation}
is an isomorphism. Indeed, in degrees $\leq n$, $\phi_n$ is identity. The same goes for $\phi_n\otimes \op{id}_M$, and hence for the induced map in homology up to degree $n-1$.

Consequently, it suffices to prove that for all $n\in\mathbb{N}$, the map $\phi_{\bullet, N_{\leq n}, M}$, is an isomorphism. We proceed by induction.

For $n=0$, the projection $\pi_0:N_{\leq 0}\to H_0(N_{\leq 0})$
is a weak equivlence. As $M$ is a cofibrant $\cA$-module, $-\otimes_\cA M$ preserves weak equivalences. Consequenlty, the composition
$$N_{\leq 0}\otimes_{\cA}M\xrightarrow{\pi_0\otimes M} H_0(N_{\leq 0})\otimes_{\cA}M\xrightarrow{\phi\circ\pi_0}  H_0(N_{\leq 0})\otimes_{H_0\cA}H_0M$$
is a weak equivalence as well, and the induced isomorphism in homology is exactly the inverse of  (\ref{propo}).

 Suppose now that $\phi_{\bullet,N_{\leq n-1}, M}$ is isomorphism for an $n\in\mathbb{N}$.  Kernel of the fibration $p_n: N_{\leq n}\twoheadrightarrow N_{\leq n-1}$ is the complex 
$$\op{Ker}p_n=(\ldots\to 0\to \op{im}(d^N_{n+1})\hookrightarrow \op{ker}(d_n^N)\to 0\to\ldots)$$
concentrated in degrees $n$ and $n+1$. There is the evident weak equivalence
$\pi_n: \op{Ker}p_n\to  H_n(N)[n]$.  As $M$ is flat and cofibrant, $$\op{Ker}p_n\otimes_{\cA} M\xrightarrow{\op{ker}(p_n)\otimes M}  N_{\leq n}\otimes_{\cA}M\xrightarrow{p_n\otimes M} N_{\leq n-1}\otimes_{\cA} M$$
is a homotopy fiber sequence. 
From the associated long exact sequence in homology, we read the following:
\begin{itemize}
\item For $0<i<n$, $H_i(p_n\otimes M)$ is an isomorphism. It now follows from the induction hypothesis that $\phi_{i,N_{\leq n}, M}$ is isomorphism. Degree $i=0$ is subtle as the long exact sequence in homology confirms only that  $H_0(p_n\otimes M)$ is injective. The map is surjective as $p_n$ is degreewise surjective, and functors $-\otimes_{\cA} M:\tt Mod(\cA)\to Mod(\cA)$, $H_0:{\tt Mod(\cA)\to Mod}(H_0\cA)$ are both left adjoint, thus right exact. Functor $-\otimes_{\cA} M$ is left adjoint as $\tt Mod(\cA)$ is a closed monoidal category, and $H_0$ is the left adjoint of the inclusion functor ${\tt Mod}(H_0\cA) \hookrightarrow {\tt Mod}(\cA)$.
\item $H_n(\op{ker}p_n\otimes M)$ is an isomorphism. Let $\psi$ denote the composed quasi-isomorphism
$$\op{Ker}p_n\otimes_{\cA}M\xrightarrow{\pi_n\otimes M}  H_n(N)[n]\otimes_{\cA}M\xrightarrow{(\phi\circ\pi_0)[n]}   H_n(N)[n]\otimes_{H_0\cA}H_0M.$$
It is straightforwardly verified that  $\phi_{n,N_{\leq n}, M}=H_n(\op{ker}p_n\otimes M)\circ H_n(\psi)^{-1},$
hence isomorphism.
\item For $i>m$, $H_i(N_{\leq n}\otimes_{\cA} M)=0$, and $\phi_{n,N_{\leq n}, M}:0\to 0$ is trivially isomorphism.
\end{itemize}

%
\end{proof}
Finally, we have shown that notions of flatness and strong flatness are equivalent. The story of flat modules summarizes as:
\begin{theo}
For $M\in\tt Mod(\cA)$, the following properties are equivalent:
\begin{enumerate}
\item $M$ is flat,
\item the derived tensor product $-\otimes_\cA^\mathbb{L} M$ preserves homotopy fiber sequences,
\item $M$ is strongly flat.
\end{enumerate}
\end{theo}

\begin{cor}\label{last_cor}
Given $M,N\in\tt Mod(\cA)$, if $M$ is a flat $\cA$-module, the above Tor spactral sequence collapses on the second page, yielding a natural isomorphism
$$(H M\otimes_{H\cA}H N)_q\cong H_q(M\otimes^\mathbb{L}_\cA N).$$
\end{cor}
\begin{proof}
Since $M$ is flat it is also strongly flat. For strongly flat modules, the above property is proven in the subsection \ref{flat1}.
\end{proof}
\section{Outlook}\label{outlook}

This text is part of a program which aims to establish homotopical algebraic geometry over differential operators as a natural framework for partial differential equations and their symmetries \cite{CG1,P11}. For this we have to show in particular that the triplet $(\tt DG\cD M, DG\cD M, DG\cD A)$ together with \'etale coverings and smooth morphisms, is not only a homotopical algebraic but even a homotopical algebraic geometric context. This includes proving that, in our specific environment, flat (resp., \'etale) morphisms are the same as strongly flat (resp., strongly \'etale). The proofs require that Quillen's Tor spectral sequence be valid in the $\cD$-geometric setting. The latter has been proved in this work, and we expect to use it to complete the proof that solid concepts of derived stack and geometric derived stack do exist in our homotopical $\cD$-geometric setting. Viewed from a broader perspective, the present work is also part of an effort to strengthen the role of the functor-of-points as a fundamental approach in derived algebraic geometry \cite{GroFun2,TV05,Paugam1} and colored supergeometry \cite{Berezinian,Z2nManifolds,LocalForms,Integration}.

\footnotesize{{}
\vfill
{\emph{email:} {\sf alisa.govzmann@uni.lu}; \emph{email:} {\sf damjan.pistalo@uni.lu}; \emph{email:} {\sf norbert.poncin@uni.lu}.}

\end{document}